\documentclass[12pt]{article}

\usepackage{amsmath, amsthm, amssymb, mathtools, extarrows}
\usepackage{graphicx}
\usepackage{hyperref}
\usepackage{cleveref}
\usepackage{multirow}
\usepackage{caption}
\usepackage{ulem}
\usepackage{enumitem}
\usepackage{booktabs}
\usepackage[numbers]{natbib}
\usepackage{authblk}
\usepackage[margin=1in]{geometry}
\usepackage{tikz}
\usetikzlibrary{positioning, decorations, shapes.geometric, arrows.meta}

\newtheorem{theorem}{Theorem}[section]
\newtheorem{lemma}[theorem]{Lemma}
\newtheorem{proposition}[theorem]{Proposition}
\newtheorem{corollary}[theorem]{Corollary}

\theoremstyle{definition}
\newtheorem{definition}[theorem]{Definition}

\newtheorem{remark}[theorem]{Remark}

\numberwithin{equation}{section}

\newcommand{\ee}{\mathrm{e}} % Natural constant
\newcommand{\im}{\mathrm{i}} % Imaginary unit
\newcommand{\RE}{\operatorname{Re}} % Real part
\newcommand{\diag}{\operatorname{diag}} % Diagonal matrix
\newcommand{\cj}{\overline} % Conjugate
\newcommand{\tp}{\mathrm{T}} % Transpose
\newcommand{\cjtp}{\mathrm{H}} % Conjugate transpose
\newcommand{\gal}{\mathrm{Gal}} % Galois
\newcommand{\Q}{\mathbb{Q}} % Rational numbers
\newcommand{\Z}{\mathbb{Z}} % Integers
\newcommand{\tmod}{\mathrm{mod}} % Modulo

\title{
    Substructure Analysis and Cycle Enumeration Methods \\
    for Oriented Graphs Based on Parameterizing \\
    Hermitian Laplacian Matrices by Galois Conjugates
}

\author[1,*]{Silin Huang}
\affil[1]{
    \textit{College of Information Science and Technology, Jinan University,}\authorcr
    \textit{Guangzhou, Guangdong, 510632, China}\authorcr
    \texttt{huangsilin@stu.jnu.edu.cn}
}

\date{}

\begin{document}

\maketitle

\begin{abstract}
This paper investigates the principal minors
of a parameterized Hermitian Laplacian matrix for oriented graphs.
Particularly, we focus on the properties of the matrix
for parameters chosen as Galois conjugates
of a primitive $p$th root of unity, where $p$ is an odd prime.
We demonstrate that under this condition,
the product of the corresponding Hermitian Laplacian determinants
is an integer power of $p$.
This algebraic property forms the basis
for a method to enumerate non-vanishing unicyclic graph components
within certain substructures.
The study is situated within a framework
where a variable unit-modulus complex parameter is introduced
into the Hermitian Laplacian matrix,
which also allows for an examination of relationships
among principal minors under different parameters.
Our analysis adopts the concept of substructures,
defined as vertex-edge pairs $(V',E')$ where edges in $E'$
are not restricted to connecting vertices within $V'$.


\end{abstract}

{
\small
\setlength{\parskip}{1pt}
\tableofcontents
}

\section{Introduction}
This paper investigates the properties of
a parameterized Hermitian Laplacian matrix for oriented graphs,
especially when the complex parameter
is chosen as a primitive $p$th root of unity,
where $p$ is an odd prime, and its Galois conjugates.
The primary result of this paper is a generalized theorem
showing that the product of Hermitian Laplacian determinants
corresponding to a set of Galois conjugate parameters
is an integer power of $p$.
This property provides a method
for enumerating non-$p$-vanishing
unicyclic graph components within certain substructures.
A unicyclic graph is non-$p$-vanishing
if its cycle is an oriented $k$-cycle with $g$ negative edges
such that $(k-2g)\bmod p\ne 0$.

Our analysis is built upon the concept of substructures,
a notion introduced by Bapat et al. [\citenum{bapat}]
for real-valued matrices,
on which combinatorial interpretation of Laplacian minors
often relies.
A substructure is defined as a vertex-edge pair $(V',E')$,
which relaxes the traditional subgraph constraint
by permitting edges in $E'$ to connect to vertices
outside of $V'$.

The study of Hermitian matrices associated with graphs has seen
considerable development.
A branch of research began with Liu and Li [\citenum{liu}],
who introduced the Hermitian adjacency matrix,
where directed edges are weighted by $\pm\im$.
Guo and Mohar [\citenum{guo}] gave a thorough introduction
to the spectra of this matrix.
Building upon this,
Yu et al. [\citenum{yu-1}, \citenum{yu-2}]
investigated the singularity and determinant expressions for
the corresponding Hermitian Laplacian, linking them to
topological properties such as cycle positivity or negativity.
Sarma [\citenum{sarma}] later extended this work
by providing expressions for minors of
these first-kind matrices and
further categorized mixed cycles into several types.

Another line of investigation was started by Mohar [\citenum{mohar}],
who proposed the idea of selecting the most appropriate
complex parameter for the Hermitian adjacency matrix
and then introduced the concept of
Hermitian adjacency matrix of the second kind,
which is parameterized by the $6$th root of unity
$\omega_6=(1+\im\sqrt{3})/2$.
This choice is motivated by one of its algebraic properties:
$\omega_6+\cj{\omega_6}=\omega_6\cj{\omega_6}=1$,
which unifies the representation of digons
and undirected edges of a mixed graph.
Recently, Xiong et al. [\citenum{sk-1}, \citenum{sk-2}]
adopted this construction and presented some expressions for
principal and non-principal minors of this matrix.
However, these results leave the properties
for other parameters as an open question.

The motivation for exploring parameters beyond the commonly used ones
becomes natural and
can be illustrated by the advantages of the $5$th roots of unity.
In the framework utilizing the $6$th root of unity,
as explored by Xiong et al. [\citenum{sk-1}, \citenum{sk-2}],
unicyclic substructures contribute integer values
$1,3,4,\text{ or }0$ to the Hermitian Laplacian principal minors.
Notably, a contribution of $1$ from a cycle is
algebraically indistinguishable from the contribution
of a rootless tree.
In contrast, setting $\omega$ to a $5$th root of unity
maps the cycle contributions to
$\alpha=\sqrt{5}\phi,\,\beta=\sqrt{5}\phi^{-1}
\in\Q(\sqrt{5})$, and $0$,
where $\phi=(1+\sqrt{5})/2$ is the golden ratio.
This not only provides a different set
of values but also ensures that
all non-vanishing cycle contributions are
separable from those of rootless trees.
This algebraic separability was the starting point of our investigation.

The paper also demonstrates the new capabilities
unlocked by treating the complex parameter as a variable.
We show that by doing so,
we can alter the determinant contributions of cycles
(which now have a simpler algebraic expression $4\sin^2[(k-2g)\pi/p]$),
allowing for more results
not readily available with fixed-parameter matrices.
These results include a general unicyclic graph enumeration method
based on Galois conjugates,
an analysis of relationships among Hermitian Laplacian principal minors
under certain different parameters,
and a new analog of the Matrix-Tree Theorem.

The remainder of this paper is organized as follows.
\begin{itemize}
    \item In Section \ref{sec:preliminaries},
    we establish the fundamental definitions,
    including the parameterized Hermitian Laplacian matrix
    and the concept of graph substructures.

    \item Section \ref{sec:main_results}
    presents our main theoretical results,
    detailing the properties
    of the parameterized Hermitian Laplacian matrix,
    which serves as the foundation for our applications.

    \item Section \ref{sec:applications}
    demonstrates the utility of this framework
    through several applications.
    This includes a unicyclic graph enumeration method
    using a $5$th root of unity and its Galois conjugate,
    followed by the generalization of this method to any odd prime $p$.
    We then explore relationships among
    Hermitian Laplacian principal minors for various roots of unity
    and provide a new analog of the Matrix-Tree Theorem.

    \item Section \ref{sec:conclusion}
    concludes the paper with a summary of our findings.
\end{itemize}

\section{Preliminaries}\label{sec:preliminaries}
In this section, we introduce the fundamental concepts
and notation used throughout the paper.
We begin by defining a family
of parameterized Hermitian matrices:
Hermitian incidence, adjacency, and Laplacian matrix,
for oriented graphs,
in Section \ref{sec:hmt-inc-and-hmt-adj} and \ref{sec:hmt-lap}.
The definition of these matrices
is essentially similar to those with fixed parameters.
We then introduce the important concept of a graph substructure
in Section \ref{ss:def-ss},
which relaxes the traditional definition of a subgraph,
and we extend the matrix definitions to these objects
in Section \ref{sec:ss-matrices}.
Finally, we recall the Cauchy-Binet formula
in Section \ref{sec:c-b-formula}
which is essential for our main theorems.

For simplicity,
although these matrices are parameterized
and differ in some ways from the original Hermitian matrices,
we will omit the term ``parameterized''
in the rest of the paper.

\begin{definition}
    The \textit{oriented graphs} we discuss in this paper
    are graphs with directed edges ($u\rightarrow v$)
    which contain no loops ($u\rightarrow u$),
    no multiple arcs ($u\xrightarrow{e_1}v,u\xrightarrow{e_2}v,e_1\ne e_2$)
    and no digons ($u\leftrightarrow v$).
\end{definition}

\subsection{Hermitian Incidence Matrix and Hermitian Adjacency Matrix}\label{sec:hmt-inc-and-hmt-adj}
\begin{definition}
    Suppose that $G(V,E)$ is an oriented graph
    of order $|V|=n$ and size $|E|=m$.
    A \textit{Hermitian incidence matrix} of $G$ is an $n\times m$ matrix
    $S_\omega(G) = (s_{ve})$ with entries
    \begin{equation}
        s_{ve}=
        \begin{cases}
            -\omega, &\exists u\in V,\, v\xrightarrow{e}u, \\
            1, &\exists u\in V,\, v\xleftarrow{e}u, \\
            0, &\text{otherwise},
        \end{cases}
    \end{equation}
    where $\omega$ is a complex number with modulus $1$.
\end{definition}

The selection of $\omega$ introduces variability,
thereby leading to multiple possible matrix descriptions for
the same graph.
Therefore, it is important to clearly specify this parameter
when referring to such matrices in the following contexts.

\begin{definition}
    A \textit{Hermitian adjacency matrix} of
    oriented graph $G$ is an $n\times n$ matrix
    $A_\omega(G) = (a_{uv})$ with entries
    \begin{equation}
    a_{uv}=
    \begin{cases}
        \omega, &u\rightarrow v, \\
        \cj\omega, &u\leftarrow v, \\
        0, &\text{otherwise},
    \end{cases}
    \end{equation}
    where $\omega$ is a complex number with modulus $1$.
\end{definition}
The construction of this matrix
was first introduced by Mohar [\citenum{mohar}].

It is worth noting that, on one hand, setting $\omega = 1$ results
in the familiar original forms of the \textit{incidence matrix} and \textit{adjacency matrix},
which are capable of representing the relationships
between vertices and edges in unoriented graphs;
on the other hand, setting $\omega = \im$ gives rise
to the complex matrices commonly referred
to as the \textit{Hermitian incidence matrix} and \textit{Hermitian adjacency matrix},
which effectively describe the relationships
between vertices and edges in oriented graphs.

\subsection{Hermitian Laplacian Matrix}\label{sec:hmt-lap}
\begin{definition}
    Let $G=(V,E)$ be an oriented graph of order $|V|=n$.
    A \textit{Hermitian Laplacian matrix} $L_\omega(G)$ of $G$ is
    an $n\times n$ matrix defined as
    $L_\omega(G)=D(G)-A_\omega(G)$, where $D(G)$ is the degree diagonal matrix of $G$,
    given explicitly by $D(G)\coloneqq\diag(\deg(u))_{u\in V}$.
    Here the \textit{degree} of a vertex $u$,
    denoted by $\deg(u)$,
    is the sum of in-degree and out-degree of $u$
    and is equal to the number of edges connected to $u$
    in the underlying graph of $G$.
\end{definition}

For better clarity, the entry-wise definition
of $L_\omega(G) = (l_{uv})$ can be expressed as follows.
\begin{equation}
    l_{uv}=
    \begin{cases}
        -\omega, &u\rightarrow v, \\
        -\cj\omega, &u\leftarrow v, \\
        \deg(u), &u=v, \\
        0, &\text{otherwise}.
    \end{cases}
\end{equation}

\begin{remark}
    Hermitian Laplacian matrix $L_\omega(G)$ is Hermitian,
    i.e. $L_\omega(G)=(L_\omega(G))^\cjtp$,
    where $(\,\cdot\,)^\cjtp$ denotes the conjugate transpose of a complex matrix.
\end{remark}

We now present an important proposition
that clarifies the relationship between $L_\omega$ and $S_\omega$.
\begin{proposition}
    Let $G$ be an oriented graph and $\omega$ be a complex number with modulus $1$.
    Then $L_\omega(G)=S_\omega(G)(S_\omega(G))^\cjtp$.
\end{proposition}
\begin{proof}
    Suppose $L_\omega(G)=(l_{uv})_{n\times n}$ and $S_\omega(G)(S_\omega(G))^\cjtp=(p_{uv})_{n\times n}$.

    For $u=v$, we have
    \begin{equation}
        p_{uv}=\sum_{e\in E}s_{ue}\cj{s_{ue}}=\sum_{e\in E}|s_{ue}|^2=
        \sum_{e\in E}\begin{cases}
            1, &e\text{ is relevant to }u, \\
            0, &\text{otherwise}
        \end{cases}=\deg(u).
    \end{equation}

    For $u\ne v$, we have
    \begin{equation}
        \begin{aligned}
            \displaystyle p_{uv}=\sum_{e\in E}s_{ue}\cj{s_{ve}}
            &=\sum_{e\in E}\begin{cases}
                -\omega\cdot\cj{1}, &u\xrightarrow{e}v, \\
                1\cdot\cj{-\omega}, &u\xleftarrow{e}v, \\
                0, &u\text{ and }v\text{ are not connected by }e
            \end{cases} \\
            &=\begin{cases}
                -\omega, &u\rightarrow v, \\
                -\cj\omega, &u\leftarrow v, \\
                0, &u\text{ and }v\text{ are not connected}.
            \end{cases}
        \end{aligned}
    \end{equation}

    Since $l_{uv}=p_{uv}$ holds for all $u,v$, we have $L_\omega(G)=S_\omega(G)(S_\omega(G))^\cjtp$.
\end{proof}

\subsection{Substructures of Oriented Graphs}\label{ss:def-ss}
\begin{definition}
    We define a \textit{substructure} of a graph $G$
    to be an object $G'(V',E')$,
    where $V'$ and $E'$ are subsets of
    the vertex set $V$ and the edge set $E$ respectively.
\end{definition}

\begin{remark}
    A substructure need not be a subgraph, since
    it is valid for $E'$ to contain edges that have endpoints
    outside of $V'$,
    which makes $(V',E')$ fail to be a graph.
\end{remark}

The concept of a substructure was first introduced
by Bapat et al. [\citenum{bapat}]
and can be seen as a generalization of the concept of a subgraph.

An example of substructures is shown in Figure \ref{fig:ss}.

\begin{figure}[htbp]
\centering
\begin{tikzpicture}[
    vertex/.style={circle, draw, thick, minimum size=0.7cm, font=\small},
    v_prime/.style={vertex, fill=blue!20, draw=blue!50!black, very thick},
    e_prime/.style={blue!70!black, very thick, -{Stealth[length=2mm, width=1.5mm]}},
    e_prime_gray/.style={gray!70!black, very thick, -{Stealth[length=2mm, width=1.5mm]}},
    v_other/.style={vertex, fill=gray!10, draw=gray!60},
    e_other/.style={gray!50, thin, -{Stealth[length=2mm, width=1.5mm]}},
    label_style/.style={font=\small}
]
% --- Graph G ---
\begin{scope}[xshift=0cm, yshift=0cm]
    \node[vertex] (g1) at (0,2) {1};
    \node[vertex] (g2) at (2,2) {2};
    \node[vertex] (g3) at (0,0) {3};
    \node[vertex] (g4) at (2,0) {4};
    \node[vertex] (g5) at (4,1) {5};
    \draw[e_prime_gray, shorten >=1pt, shorten <=1pt] (g1) to (g2);
    \draw[e_prime_gray, shorten >=1pt, shorten <=1pt] (g1) to (g3);
    \draw[e_prime_gray, shorten >=1pt, shorten <=1pt] (g2) to (g4);
    \draw[e_prime_gray, shorten >=1pt, shorten <=1pt] (g3) to (g4);
    \draw[e_prime_gray, shorten >=1pt, shorten <=1pt] (g4) to (g5);
    \draw[e_prime_gray, shorten >=1pt, shorten <=1pt] (g2) to (g5);
    \node[label_style, above=0.2cm of g2, align=center] {Graph $G$};
\end{scope}
% --- Substructure 1 ---
\begin{scope}[xshift=6cm, yshift=0cm]
    \node[v_prime] (s1_1) at (0,2) {1};
    \node[v_other] (s1_2) at (2,2) {2};
    \node[v_prime] (s1_3) at (0,0) {3};
    \node[v_other] (s1_4) at (2,0) {4};
    \node[v_other] (s1_5) at (4,1) {5};
    \draw[e_other, shorten >=1pt, shorten <=1pt] (s1_1) to (s1_2);
    \draw[e_prime, shorten >=1pt, shorten <=1pt] (s1_1) to (s1_3);
    \draw[e_other, shorten >=1pt, shorten <=1pt] (s1_2) to (s1_4);
    \draw[e_other, shorten >=1pt, shorten <=1pt] (s1_3) to (s1_4);
    \draw[e_other, shorten >=1pt, shorten <=1pt] (s1_4) to (s1_5);
    \draw[e_other, shorten >=1pt, shorten <=1pt] (s1_2) to (s1_5);
    \node[label_style, above=0.2cm of s1_2, align=center] {Substructure $G'_1$ (a subgraph)};
    \node[below=0.1cm of s1_4, align=center, font=\tiny] {$V'=\{1,3\}$ \\ $E'=\{1\to 3\}$};
\end{scope}
% --- Substructure 2 ---
\begin{scope}[xshift=0cm, yshift=-5cm]
    \node[v_other] (s2_1) at (0,2) {1};
    \node[v_prime] (s2_2) at (2,2) {2};
    \node[v_other] (s2_3) at (0,0) {3};
    \node[v_prime] (s2_4) at (2,0) {4};
    \node[v_other] (s2_5) at (4,1) {5};
    \draw[e_other, shorten >=1pt, shorten <=1pt] (s2_1) to (s2_2);
    \draw[e_other, shorten >=1pt, shorten <=1pt] (s2_1) to (s2_3);
    \draw[e_prime, shorten >=1pt, shorten <=1pt] (s2_2) to (s2_4);
    \draw[e_other, shorten >=1pt, shorten <=1pt] (s2_3) to (s2_4);
    \draw[e_prime, shorten >=1pt, shorten <=1pt] (s2_4) to (s2_5);
    \draw[e_other, shorten >=1pt, shorten <=1pt] (s2_2) to (s2_5);
    \node[label_style, above=0.2cm of s2_2, align=center] {Substructure $G'_2$};
    \node[below=0.1cm of s2_4, align=center, font=\tiny] {$V'=\{2,4\}$ \\ $E'=\{2\to 4, 4\to 5\}$};
\end{scope}
% --- Substructure 3 ---
\begin{scope}[xshift=6cm, yshift=-5cm]
    \node[v_prime] (s3_1) at (0,2) {1};
    \node[v_other] (s3_2) at (2,2) {2};
    \node[v_other] (s3_3) at (0,0) {3};
    \node[v_other] (s3_4) at (2,0) {4};
    \node[v_prime] (s3_5) at (4,1) {5};
    \draw[e_prime, shorten >=1pt, shorten <=1pt] (s3_1) to (s3_2);
    \draw[e_other, shorten >=1pt, shorten <=1pt] (s3_1) to (s3_3);
    \draw[e_other, shorten >=1pt, shorten <=1pt] (s3_2) to (s3_4);
    \draw[e_other, shorten >=1pt, shorten <=1pt] (s3_3) to (s3_4);
    \draw[e_prime, shorten >=1pt, shorten <=1pt] (s3_4) to (s3_5);
    \draw[e_other, shorten >=1pt, shorten <=1pt] (s3_2) to (s3_5);
    \node[label_style, above=0.2cm of s3_2, align=center] {Substructure $G'_3$};
    \node[below=0.1cm of s3_4, align=center, font=\tiny] {$V'=\{1,5\}$ \\ $E'=\{1\to 2, 5\to 4\}$};
\end{scope}
\end{tikzpicture}
\caption{Examples of substructures:
a subgraph, a substructure with an edge outside the vertex set
and a substructure with two edges outside the vertex set.
}
\label{fig:ss}
\end{figure}

\subsubsection{Substructure Connectedness}
For a substructure, we define its \textit{connectedness}.
Although this property has a well-established and generalized
definition for subgraphs,
it requires redefinition for substructures.

\begin{definition}
    For a substructure $G'(V',E')$ of a graph $G(V, E)$:
    \begin{itemize}
    \item For two vertices $u, v$ in $V'$,
    if there exists vertices $u_1,\ldots,u_n\in V'$
    and edges $e_1,\ldots,e_{n-1}\in E'$
    which form a path
    $u=u_1\overset{e_1}{\sim}u_2\overset{e_2}{\sim}\cdots\overset{e_{n-1}}{\sim}u_n=v$
    connecting them,
    then $u$ and $v$ are said to be \textit{connected}.
    Here $\sim$ represents to a directed edge,
    ignoring the direction.

    \item If any two vertices in $V'$ are connected,
    then substructure $G'(V',E')$ is defined to be \textit{connected}.
    
    \item If $V'$ contains only one vertex, it is also defined as connected.
    If $V'$ contains no vertices, it is defined to be not connected.
    \end{itemize}
\end{definition}

\begin{definition}
    According to the connectivity between vertices,
    a substructure can be partitioned into parts,
    in which every two vertices are connected,
    and vertices from different parts are not connected.
    We call those parts
    \textit{connected components} (in short \textit{components}).
\end{definition}

\subsubsection{Substructure Regularness}
Building on Bapat et al. [\citenum{bapat}],
who defined \textit{nonsingular substructures}
as substructures that have square incidence matrices
and have nonzero Laplacian determinants,
we define the \textit{regularness} and \textit{all-regularness}
of a substructure.

\begin{definition}
    A \textit{regular substructure} is a substructure $G'(V',E')$
    such that $|V'|=|E'|$.
\end{definition}

\begin{definition}\label{def:all-regular-ss}
    An \textit{all-regular substructure} is a substructure $G'(V',E')$
    such that its every connected component $(V'_i,E'_i)\subseteq (V',E')$
    satisfies $|V'_i|=|E'_i|$.
\end{definition}

\subsection{Substructure Matrices}\label{sec:ss-matrices}
For the concept of a substructure,
we can similarly define its Hermitian incidence matrix
and Hermitian Laplacian matrix as follows.

\begin{definition}
    Hermitian Laplacian matrix of a substructure $G'(V',E')\subseteq G(V,E)$,
    written $L_\omega(G')=(l_{uv})$, is a $|V'|\times|V'|$ matrix
    defined as:
    
    for $u=v$,
    \begin{equation}
        l_{uv}=\deg_{G'}(u)\coloneqq\#\{e\in E':u\xrightarrow{e}w\text{ or }u\xleftarrow{e}w\};
    \end{equation}

    for $u\ne v$,
    \begin{equation}
        l_{uv}=\begin{cases}
            -\omega, &\exists e\in E', u\xrightarrow{e} v, \\
            -\cj\omega, &\exists e\in E', u\xleftarrow{e} v, \\
            0, &\text{otherwise}.
        \end{cases}
    \end{equation}
\end{definition}

\begin{definition}
    Hermitian incidence matrix of a substructure $G'(V',E')\subseteq G(V,E)$
    is defined to be the submatrix formed by
    the rows corresponding to $V'$
    and the columns corresponding to $E'$
    in $S_\omega(G(V,E))$.
    
    For better clarity, the entry-wise definition
    of $S_\omega(G')=(l_{ve})$ can be expressed as follows.
    \begin{equation}
        s_{ve}=\begin{cases}
            -\omega, &\exists u\in V,\, v\xrightarrow{e}u, \\
            1, &\exists u\in V,\, v\xleftarrow{e}u, \\
            0, &e\text{ is not relevant to }v.
        \end{cases}
    \end{equation}
\end{definition}

\subsection{Principal Minors and Cauchy-Binet Formula}\label{sec:c-b-formula}
\textit{Principal minors} of a matrix
are the determinants of submatrices
formed by rows and columns of the same indices from the matrix
and can be decomposed using the Cauchy-Binet formula.

\begin{theorem}[Cauchy-Binet formula, {[\citenum{c-b-formula}]}]
    Let $\mathbb{F}$ be a number field.
    Let $A\in\mathbb{F}^{n\times m}$ and $B \in \mathbb{F}^{m\times n}$ be two matrices on $\mathbb{F}$.
    Let $S\subseteq\{1,2,\ldots,m\}$ be a subset of size $n$.
    Denote by $A[-;S]$ the $n\times n$ submatrix of $A$ consisting of
    every original rows but only the columns indexed by $S$,
    and by $B[S;-]$ the $n\times n$ submatrix of $B$ consisting of
    every original columns but only the rows indexed by $S$.
    Then the determinant of the matrix product $AB$ is given by
    \begin{equation}
        \det(AB)=\sum_{\substack{S\subseteq\{1,2,\ldots,m\}\\|S|=n}}\det(A[-;S])\det(B[S;-]).
    \end{equation}
\end{theorem}

In the remainder of this paper, we will likewise use the notation $A[R;C]$
to denote a submatrix of $A$ consisting of rows indexed by $R$ and columns indexed by $C$.
If $R=C$ (which makes the determinant of this submatrix a principal minor),
we may write $A[R]$ or $A[C]$ for convenience.

\section{Main Results}\label{sec:main_results}
\subsection{Hermitian Laplacian Matrix Properties on Substructures}
\subsubsection{Substructure Matrix Properties}
The following lemma shows that
the concept of Hermitian Laplacian matrix of graph substructures
is suitable for representing them,
since one of the most important properties of the matrix
still holds.

\begin{lemma}\label{lem:lapl_decomp_ss}
    Let $G$ be an oriented graph, $G'(V',E')$ be its substructure and $\omega$ be a complex number with modulus $1$.
    Then $L_\omega(G')=S_\omega(G')(S_\omega(G'))^\cjtp$.
\end{lemma}
\begin{proof}
    Suppose $L_\omega(G')=(l_{uv})_{n\times n}$, $S_\omega(G')=(s_{uv})_{n\times m}$
    and $S_\omega(G')(S_\omega(G'))^\cjtp=(p_{uv})_{n\times n}$.

    For $u=v$, we have
    \begin{equation}
        \begin{aligned}
            p_{uv}&=\sum_{e\in E'}s_{ue}\cj{s_{ue}}=\sum_{e\in E'}|s_{ue}|^2 \\
            &=\sum_{e\in E'}\begin{cases}
                1, &e\text{ is relevant to }u\text{ and }u\in V', \\
                0, &\text{otherwise}
            \end{cases}=\deg_{G'}(u).
        \end{aligned}
    \end{equation}

    For $u\ne v$, we have
    \begin{equation}
        \begin{aligned}
            \displaystyle p_{uv}=\sum_{e\in E'}s_{ue}\cj{s_{ve}}
            &=\sum_{e\in E'}\begin{cases}
                -\omega\cdot\cj{1}, &u\xrightarrow{e}v\text{ and }u,v\in V', \\
                1\cdot\cj{-\omega}, &u\xleftarrow{e}v\text{ and }u,v\in V', \\
                0\cdot\cj\ast\text{ or }\ast\cdot\cj 0=0, &\text{an endpoint of }e\text{ is not in }V',\\
                0\cdot\cj 0=0, &u\text{ and }v\text{ are not connected by }e
            \end{cases} \\
            &=\begin{cases}
                -\omega, &u\rightarrow v\text{ and }u,v\in V', \\
                -\cj\omega, &u\leftarrow v\text{ and }u,v\in V', \\
                0, &u\notin V'\text{ or }v\notin V'\\&\text{or }u,v\in V'\text{ but }u\text{ and }v\text{ are not connected}.
            \end{cases}
        \end{aligned}
    \end{equation}

    It is obvious that $l_{uv}=p_{uv}$ holds for all $u,v$. Thus $L_\omega(G')=S_\omega(G')(S_\omega(G'))^\cjtp$.
\end{proof}

The following lemma presents a useful property of this matrix:
its determinant remains unchanged
upon the addition of pendant edges.

\begin{lemma}\label{lem:pendant_doesnt_change_detL}
    For a substructure, adding a pendant edge
    with a pendant vertex attached to it
    does not change the determinant of its Hermitian Laplacian matrix.
\end{lemma}
\begin{proof}
    Suppose we have a substructure $G'(V',E')$ with $|V'|=n$
    whose Hermitian Laplacian matrix is denoted by $L$.
    Then the addition of a pendant edge $e=v_iv_{n+1}$ where
    a new pendant vertex $v_{n+1}$ connects to the vertex $v_i$
    causes the following changes to happen:
    \begin{equation}
        L=\left[ \begin{matrix}
            d_1&		l_{12}&		\cdots&		l_{1i}&		\cdots&		l_{1n}\\
            l_{21}&		d_2&		\cdots&		l_{2i}&		\cdots&		l_{2n}\\
            \vdots&		\vdots&		&		\vdots&		&		\vdots\\
            l_{i1}&		l_{i2}&		\cdots&		d_i&		\cdots&		l_{in}\\
            \vdots&		\vdots&		&		\vdots&		&		\vdots\\
            l_{n1}&		l_{n2}&		\cdots&		l_{ni}&		\cdots&		l_{nn}\\
        \end{matrix} \right]
        \rightarrow
        L_{\text{new}}=
        \left[ \begin{matrix}
            d_1&		l_{12}&		\cdots&		l_{1i}&		\cdots&		l_{1n}&		0\\
            l_{21}&		d_2&		\cdots&		l_{2i}&		\cdots&		l_{2n}&		0\\
            \vdots&		\vdots&		&		\vdots&		&		\vdots&		\vdots\\
            l_{i1}&		l_{i2}&		\cdots&		d_i+1&		\cdots&		l_{in}&		l_{i,n+1}\\
            \vdots&		\vdots&		&		\vdots&		&		\vdots&		\vdots\\
            l_{n1}&		l_{n2}&		\cdots&		l_{ni}&		\cdots&		l_{nn}&		0\\
            0&		0&		\cdots&		\overline{l_{i,n+1}}&		\cdots&		0&		1\\
        \end{matrix} \right],
    \end{equation}
    where $l_{i,n+1}\in\{-\omega,-\cj\omega\}$.
    By adding $-l_{i,n+1}$ times the $(n+1)$-th row
    to the $i$-th row of $L_{\text{new}}$, we obtain
    \begin{equation}
        L_{\text{new}}\rightarrow\left[ \begin{matrix}
            d_1&		l_{12}&		\cdots&		l_{1i}&		\cdots&		l_{1n}&		0\\
            l_{21}&		d_2&		\cdots&		l_{2i}&		\cdots&		l_{2n}&		0\\
            \vdots&		\vdots&		&		\vdots&		&		\vdots&		\vdots\\
            l_{i1}&		l_{i2}&		\cdots&		d_i&		\cdots&		l_{in}&		0\\
            \vdots&		\vdots&		&		\vdots&		&		\vdots&		\vdots\\
            l_{n1}&		l_{n2}&		\cdots&		l_{ni}&		\cdots&		l_{nn}&		0\\
            0&		0&		\cdots&		\overline{l_{i,n+1}}&		\cdots&		0&		1\\
        \end{matrix} \right].
    \end{equation}
    Expanding the determinant of this matrix along its $(n+1)$-th column
    immediately yields $\det L_{\text{new}}=\det L$.
\end{proof}

Applying the above lemma in reverse yields the following theorem.
\begin{theorem}
    For a substructure, adding or removing a pendant edge
    and a pendant vertex attached to it
    does not change the determinant of its Hermitian Laplacian matrix.
\end{theorem}

Applying the above theorem continuously yields the following corollary.
\begin{corollary}\label{cor:pendants_dont_change_detL}
    For a substructure, adding or removing
    multiple pendant edges and pendant vertices onto it
    does not change the determinant of its Hermitian Laplacian matrix.
\end{corollary}

\subsubsection{Determinant Factorization by Connected Components}
\begin{lemma}\label{lem:block-diag}
The Hermitian Laplacian determinant of a substructure
is equal to the product of determinants of its connected components.
Formally, for a substructure $G'=\bigcup_{i=1}^m G'_i$
with disjoint connected components $G'_i\,(i=1,2,\ldots,m)$,
we have
\begin{equation}
    \det L_\omega(G')=\prod_{i=1}^m\det L_\omega(G'_i).    
\end{equation}
\end{lemma}

\begin{proof}
    This follows directly from the block-diagonal structure of
    $L_\omega(G')$ when $G'$ is partitioned into connected components,
    where each block corresponds to a connected component.
\end{proof}

\subsubsection{Hermitian Laplacian Determinant of Rootless Trees and Normal Trees}
\begin{definition}
    A substructure $G'(V',E')$ is called a \textit{tree}
    if it is connected and contains no cycles.
    In this definition, there are two types of trees:
    \begin{enumerate}[label=(\arabic*)]
        \item A tree graph. This is what we commonly refer to as a ``tree.''
        It satisfies $|V'|=|E'|+1$.

        To avoid ambiguity, we will refer to this type of tree
        as a \textit{normal tree} in the following text.
        \item A non-graph tree substructure. It is formed by removing a pendant vertex
        (but keeping the pendant edge connected to it) from a tree graph.
        It satisfies $|V'|=|E'|$.

        To avoid ambiguity, we will refer to this type of tree
        as a \textit{rootless tree} in the following text.

        An example of these two types of trees is shown in Figure \ref{fig:trees}.
    \end{enumerate}

    \begin{figure}[htbp]
    \centering
    \begin{tikzpicture}[
        vertex/.style={circle, draw, thick, minimum size=0.7cm, font=\small},
        v_prime/.style={vertex, fill=blue!20, draw=blue!50!black, very thick},
        e_prime/.style={blue!70!black, very thick, -{Stealth[length=2mm, width=1.5mm]}},
        v_normal/.style={vertex, fill=green!20, draw=green!50!black, very thick},
        e_normal/.style={green!70!black, very thick, -{Stealth[length=2mm, width=1.5mm]}},
        v_removed/.style={vertex, fill=gray!10, draw=gray!60},
        e_removed/.style={gray!50, very thick, -{Stealth[length=2mm, width=1.5mm]}},
        label_style/.style={font=\small}
    ]
    % --- Normal Tree ---
    \begin{scope}[xshift=0cm, yshift=0cm]
        \node[v_normal] (n1) at (0,1.5) {1};
        \node[v_normal] (n2) at (-1,0) {2};
        \node[v_normal] (n3) at (1,0) {3};
        \node[v_normal] (n4) at (1,-1.5) {4};
        \draw[e_normal, shorten >=1pt, shorten <=1pt] (n1) to (n2);
        \draw[e_normal, shorten >=1pt, shorten <=1pt] (n1) to (n3);
        \draw[e_normal, shorten >=1pt, shorten <=1pt] (n3) to (n4);
        \node[label_style, above=0.3cm of n1, align=center] {A $4$-normal tree $T_1$};
    \end{scope}
    % --- Rootless Tree ---
    \begin{scope}[xshift=6cm, yshift=0cm]
        \node[v_prime] (r11) at (0,1.5) {1};
        \node[v_removed] (r12) at (-1,0) {2};
        \node[v_prime] (r13) at (1,0) {3};
        \node[v_prime] (r14) at (1,-1.5) {4};
        \draw[e_prime, shorten >=1pt, shorten <=1pt] (r11) to (r12);
        \draw[e_prime, shorten >=1pt, shorten <=1pt] (r11) to (r13);
        \draw[e_prime, shorten >=1pt, shorten <=1pt] (r13) to (r14);
        \node[label_style, above=0.3cm of r11, align=center] {A $3$-rootless tree $T'_2$};
    \end{scope}

    \begin{scope}[xshift=-1.4cm, yshift=-4cm]
        \node[v_normal] (n21) at (0,0) {1};
        \node[v_normal] (n22) at (1.5,0) {2};
        \node[v_removed] (n23) at (3,0) {3};
        \draw[e_normal, shorten >=1pt, shorten <=1pt] (n21) to (n22);
        \draw[e_removed, shorten >=1pt, shorten <=1pt] (n22) to (n23);
        \node[label_style, above=0.3cm of n22, align=center] {A 2-normal tree $T'_3$};
    \end{scope}

    \begin{scope}[xshift=4.6cm, yshift=-4cm]
        \node[v_prime] (r21) at (0,0) {1};
        \node[v_prime] (r22) at (1.5,0) {2};
        \node[v_removed] (r23) at (3,0) {3};
        \draw[e_prime, shorten >=1pt, shorten <=1pt] (r21) to (r22);
        \draw[e_prime, shorten >=1pt, shorten <=1pt] (r22) to (r23);
        \node[label_style, above=0.3cm of r22, align=center] {A 2-rootless tree $T'_3$};
    \end{scope}
    \end{tikzpicture}
    \caption{Examples of normal trees (left) and rootless trees (right).}
    \label{fig:trees}
    \end{figure}
\end{definition}

The concept of a ``rootless tree'' is
a rigorous definition of the structure
formed by removing a vertex from a classical tree
but keeping the edge attached to it.
This structure is built on Bapat et al.'s
[\citenum{bapat}] notion of a ``rootless spanning tree,''
and is fundamental in this framework.
As we will show in Theorem \ref{thm:tree-det},
it is the simplest regular substructure,
as its Hermitian Laplacian determinant is always $1$.

\begin{lemma}\label{lem:lapl_sngl_ss_regl}
    Let $G'(V',E')$ be a graph substructure and $|V'|\geqslant|E'|$.
    If the Hermitian Laplacian matrix of $G'$ is not singular,
    then every connected component of $G'$ is regular,
    i.e. $G'$ is all-regular.
\end{lemma}
\begin{proof}
    Suppose not, then there exists a component of $G'$,
    say $G'_0(V'_0,E'_0)$, that is connected but not regular.
    
    Without loss of generality, we assume $|V'_0|>|E'_0|$.
    We notice that $|V'_0|>|E'_0|+1$ is not possible,
    since connecting $|V'_0|$ vertices requires at least $|V'_0|-1$ edges.
    So the only possibility is $|V'_0|=|E'_0|+1$.
    This indicates that $G'_0$ is a normal tree.

    Obviously, a tree can be formed by continuously
    adding pendant edges and pendant vertices onto a 2-normal tree
    whose Hermitian Laplacian matrix is
    \begin{equation}
        \left[\begin{matrix}
            1 & l_{12} \\
            \cj{l_{12}} & 1
        \end{matrix}\right],
    \end{equation}
    where $l_{12}\in\{-\omega,-\cj\omega\}$.
    This matrix has determinant $0$.
    By Corollary \ref{cor:pendants_dont_change_detL}
    the determinant of the Hermitian Laplacian matrix
    of the tree is also $0$.
    By Lemma \ref{lem:block-diag},
    the determinant of $G'$ is $0$.
    This completes the proof.
\end{proof}

\begin{theorem}\label{thm:tree-det}
    If a substructure $G'(V',E')$ is a tree, then the determinant of
    its Hermitian Laplacian matrix is $1$ if it is regular;
    or $0$ if it is not regular (it is normal).
\end{theorem}
\begin{proof}
    We consider the two cases based on the substructure's regularity.

    \paragraph{Case 1: $G'$ is regular:} Then this substructure is a rootless tree.
    A rootless tree can be formed by continuously adding pendant edges and pendant vertices
    onto a 2-rootless tree whose Hermitian Laplacian matrix is
    \begin{equation}
        \left[\begin{matrix}
            1 & l_{12} \\
            \cj{l_{12}} & 2
        \end{matrix}\right],
    \end{equation}
    where $l_{12}\in\{-\omega,-\cj\omega\}$.
    This matrix has determinant $1$.
    By Corollary \ref{cor:pendants_dont_change_detL},
    the determinant of the Hermitian Laplacian matrix of the tree is also $1$.
    
    \paragraph{Case 2: $G'$ is not regular:} Then the singularness of its Hermitian Laplacian matrix
    is clear by Lemma \ref{lem:lapl_sngl_ss_regl}.
\end{proof}

\subsection{Classification of Connected Regular Substructures}
\begin{theorem}\label{thm:con-reg-ss-classification}
    A connected regular substructure is either a
    rootless tree or a unicyclic graph.
\end{theorem}
\begin{proof}
    Suppose we have a substructure $G'(V',E')$ with $|V'|=|E'|$.

    \paragraph{Case 1: $G'$ is a subgraph.} Then it is well known that
    $G'$ is a unicyclic graph.

    \paragraph{Case 2: $G'$ is not a subgraph.} Then it includes some edges
    that have one endpoint not in $V'$.
    Suppose $e$ is one of them,
    then $(V',E'\setminus\{e\})$ is also connected
    and satisfies $|V'|=|E'\setminus\{e\}|+1$,
    which indicates that it is a normal tree.
    Therefore $G'$ is a rootless tree.
\end{proof}

\begin{theorem}\label{thm:uc-det-equals-to-cyc-det}
    The determinant of the Hermitian Laplacian matrix
    of a unicyclic graph is equal to
    that of the unique cycle in it.
\end{theorem}
\begin{proof}
    The proof is straightforward and can be established by
    applying Corollary \ref{cor:pendants_dont_change_detL} on a cycle.
\end{proof}

Having established that connected regular substructures are
rootless trees and unicyclic graphs,
we now extend this classification to all-regular substructures,
which are composed of one or more such connected components.

\subsection{Rootless Forests and Classification of All-Regular Substructures}\label{ss:ars-classification}
Similar to the concept of a \textit{forest} in graph theory,
we introduce the new concept of a \textit{rootless forest}.

\begin{definition}
    If every connected component of a substructure is a rootless tree,
    then we call it a \textit{rootless forest}.
\end{definition}

By this definition,
a substructure is composed of several unicyclic subgraphs
and at most one rootless forest.

\begin{lemma}\label{lem:forest-det}
    If $G'(V',E')$ is a rootless forest, then
    \begin{equation}
        \det L_\omega(G'(V',E'))=1.
    \end{equation}
\end{lemma}
\begin{proof}
    The proof is straightforward from Theorem \ref{thm:tree-det}
    and Lemma \ref{lem:block-diag}.
\end{proof}

This implies that the determinant of any all-regular substructure
is solely determined by its unicyclic components.
Specifically, an all-regular substructure $G'$
can be uniquely decomposed into a set of unicyclic graph components
$\{C_1,\ldots,C_m\}$ and a rootless forest $F$.
By Lemma \ref{lem:block-diag}, its determinant is
\begin{equation}
    \det L_\omega(G')
    =\left(\prod_{i=1}^{m}\det L_\omega(C_i)\right)
    \cdot\underbrace{\det L_\omega(F)}_{=1}
    =\prod_{i=1}^{m}\det L_\omega(C_i).
\end{equation}

Furthermore, by Theorem \ref{thm:uc-det-equals-to-cyc-det},
the determinant of each unicyclic graph $C_i$
is equal to that of its unique cycle.
Therefore, the determinant of an all-regular substructure
only depends on the properties of cycles
in its unicyclic graph components.

\subsection{Hermitian Laplacian Matrix Principal Minor Decomposition}
The importance of the substructure lies in the fact that when applying
the Cauchy-Binet formula to a Hermitian Laplacian matrix,
we extract some submatrices of the Hermitian incidence matrix,
and the structures of the graph reflected by those submatrices
are precisely regular substructures.
Furthermore, these Hermitian incidence matrices
can convert back into the Hermitian Laplacian matrix
for graph substructures which we defined earlier.

\begin{theorem}[Hermitian Laplacian Matrix Principal Minor Decomposition Theorem]\label{thm:lapl-decomp}
    Let $G(V,E)$ be an oriented graph and
    $L_\omega$ be its Hermitian Laplacian matrix.
    For a substructure $G'(V',E')\subseteq G(V,E)$,
    denote its Hermitian Laplacian matrix by $L_\omega(G'(V',E'))$.
    Then we have
    \begin{equation}
        \displaystyle\det L_\omega[V']=\sum_{\substack{E'\subseteq E\\|E'|=|V'|}}\det L_\omega(G'(V',E')).
    \end{equation}
\end{theorem}

\begin{proof}
    Noticing that
    \begin{equation}
        L_\omega[V']=S_\omega[V';-]S_\omega[V';-]^\cjtp
    \end{equation}
    and by applying the Cauchy-Binet formula, we have
    \begin{equation}
        \begin{aligned}
            \det L_\omega[V']&=\det\left(S_\omega[V';-]S_\omega[V';-]^\cjtp\right) &\quad\\
            &=\sum_{\substack{E'\subseteq E\\|E'|=|V'|}}\det(S_\omega[V';E'])\det\left(S_\omega[V';E']^\cjtp\right) &\quad\text{(By Cauchy-Binet formula)}\\
            &=\sum_{\substack{E'\subseteq E\\|E'|=|V'|}}\left|\det\left(S_\omega[V';E']\right)\right|^2 &\quad\text{(By }\det(A^\cjtp)=(\det A)^\cjtp\text{)}\\
            &=\sum_{\substack{E'\subseteq E\\|E'|=|V'|}}\det L_\omega(G'(V',E')) &\quad\text{(By Lemma \ref{lem:lapl_decomp_ss})} \\
            &=\sum_{\substack{E'\subseteq E\\(V',E')\text{ all-regular}}}\det L_\omega(G'(V',E')). &\quad\text{(By Lemma \ref{lem:lapl_sngl_ss_regl})}
        \end{aligned}
    \end{equation}

    This completes the proof.
\end{proof}

\subsection{Walk Weight and Hermitian Laplacian Determinant of Cycles}
Similar to Xiong et al. [\citenum{sk-1}], we adopt
the concept of the \textit{weight} of a walk in an oriented graph,
which makes the analysis of
the Hermitian Laplacian determinants of cycles easier.

\begin{definition}
    A \textit{walk} in an oriented graph is defined as
    a sequence of vertices $(v_1,v_2,\ldots,v_n)$,
    where there exists a directed edge between any two consecutive vertices.
    The \textit{weight} of the walk is defined as
    $w_\omega=a_{v_1 v_2}a_{v_2 v_3}\cdots a_{v_{n-1} v_n}$,
    where $a_{ij}$ is the element at position $(i,j)$
    of the Hermitian adjacency matrix
    (or Hermitian Laplacian matrix,
    their entries of those are same here)
    of the graph.
\end{definition}

From this definition, it can be seen that the weight of a cycle
is a product of some $\omega$ and $\cj\omega$,
where whether each $\omega$ is conjugated or not depends on
the direction of the edges forming the cycle. After designating
one (``clockwise'' or ``counterclockwise'') as the positive direction,
corresponding to $\omega$;
the other becomes the negative direction,
corresponding to $\cj\omega$.

Noticing that $\cj\omega=\omega^{-1}$, the weight of a cycle
can only come from the set
$\Omega_\omega\coloneqq\{\omega^n:n\in\Z\}$.
To simplify the situation, we may choose an appropriate $\omega$
such that $\Omega_\omega$ becomes a finite set,
which is true when $(\Omega_\omega,\cdot)$ forms a cyclic group,
where $\cdot$ denotes complex multiplication.
This requires $\omega$ to be a \textit{root of unity},
specifically a \textit{primitive root of unity}, i.e.
$\omega=\ee^{\im 2\pi/p},p=a/b\,(a,b\in\Z_+,\,\gcd(a,b)=1)$.

Suppose we have already designated the positive direction
of a cycle to be the most frequently appeared direction
in that cycle.
Then, for a $k$-cycle, it must contain
$g$ negative-directed edges
and $k-g$ positive-directed edges, where
$g\in\Z\cap[0,\lfloor k/2 \rfloor]$.
Under these circumstances, the weight of this cycle is
\begin{equation}
    w_\omega=\omega^{k-g}\cj\omega^g=\omega^{k-g}\omega^{-g}=\omega^{k-2g}.
\end{equation}

\begin{definition}
    For convenience, we call a negative-directed edge
    a \textit{negative edge}.
\end{definition}

The following proposition makes calculating
the Hermitian Laplacian determinant of a cycle substructure
much easier.

\begin{proposition}\label{prp:cycle_lap}
    Let $C$ be a directed cycle. Then
\begin{equation}
    \det(L_\omega(C)) = 2-2\RE(w_\omega(C)).
\end{equation}
\end{proposition}
\begin{proof}
    The proof of this proposition
    is similar to that of Theorem 2.3 in [\citenum{sk-1}].
\end{proof}

This simple formula is central to this framework.
It transforms the problem of calculating a matrix determinant
into a simpler task of evaluating the real part
of a cycle weight, which is a product of
$\omega$ and $\cj\omega$ terms
determined by the cycle's edge orientations.

This also implies that the choice of positive direction
is inconsequential, since
\begin{equation}
    w_{\cj\omega}(C)=\cj{w_\omega(C)}\quad\Rightarrow\quad
    2-2\RE(w_{\cj\omega}(C))=2-2\RE(w_\omega(C)).
\end{equation}

\subsection{Hermitian Laplacian Determinant of Unicyclic Graphs and Its Properties}
Since it's unicyclic graph components that are dominating
the Hermitian Laplacian determinants of substructures,
we now focus on Hermitian Laplacian determinants of unicyclic graphs.

\subsubsection{A Simplified Expression for Unicyclic Graph Determinants}
In the rest of the paper,
we denote the determinant of the Hermitian Laplacian matrix
of a unicyclic graph with a $k$-cycle and $g$ negative edges by
$\ell_\omega(k,g)=2-2\RE(\omega^{k-2g})$.
The following theorem gives a simpler expression
for $\ell_\omega(k,g)$.

\begin{theorem}
    Suppose $p\in\Q_+$
    and $\omega=\ee^{\im 2\pi/p}$ is a root of unity.
    Then we have
    \begin{equation}
        \ell_\omega(k,g)=4\sin^2\left(\dfrac{k-2g}{p}\pi\right).
    \end{equation}
\end{theorem}
\begin{proof}
    Simply by Euler's formula, we have
    \begin{equation}
        \begin{aligned}
            \ell_\omega(k,g)&=2-2\RE(\omega^{k-2g}) \\
            &=2-2\RE\left(\ee^{\im\frac{2\pi}{p}(k-2g)}\right) \\
            &=2-2\cos\left[\dfrac{2\pi}{p}(k-2g)\right] \\
            &=4\sin^2\left(\dfrac{k-2g}{p}\pi\right).
        \end{aligned}
    \end{equation}
\end{proof}

\begin{proposition}
    For all $k\in\mathbb{N},g\in\Z\cap[0,\lfloor k/2 \rfloor]$ we have
    \begin{equation}
        \ell_\omega(k,g)\in[0,4].
    \end{equation}
\end{proposition}

\subsubsection{Singularity of Hermitian Laplacian Matrix of Unicyclic Graphs}\label{sss:vanish}
Yu et al. [\citenum{yu-1}, \citenum{yu-2}] characterized
the singularity of Hermitian Laplacian matrix of the first-kind
based on whether some substructures are positive.
The parameterized $\omega$ framework may allow for a clearer
understanding of when unicyclic graph substructures will
correspond to singular or non-singular Hermitian Laplacians.

\begin{theorem}
    $\ell_\omega(k,g)$ vanishes, i.e.
    \begin{equation}
        \ell_\omega(k,g)=4\sin^2\left(\dfrac{k-2g}{p}\pi\right)=0,
    \end{equation}
    when and only when
    \begin{equation}
        \dfrac{k-2g}{p}\in\Z.
    \end{equation}
    As the result, for all $g\in\Z\cap[0,\lfloor k/2 \rfloor]$,
    $\ell_\omega(k,g)$ don't vanish, when and only when
    \begin{equation}
        \dfrac{k-2g}{p}\notin\Z\quad(\forall g\in\Z\cap[0,\lfloor k/2 \rfloor])
    \end{equation}
    which holds if and only if
    \begin{enumerate}[label=(\arabic*)]
        \item $k$ is odd and $p$ is large enough so that
        $0<k-2g<p\Rightarrow 0<\frac{k-2g}{p}<1$ always holds; or
        \item $k$ is odd while $p$ is even.
    \end{enumerate}
\end{theorem}
\begin{proof}
    The proof is trivial.
\end{proof}

It's worth mentioning that if $k$ is even, then it will always be
possible for $k-2g=0$ to occur (when $g=k/2$).
In this case, the possibility of $(k-2g)/p=0\in\Z$ still exists.
In fact, an even $k$-cycle with $k/2$ negative edges
has a singular Hermitian Laplacian matrix no matter the choice of $\omega$.
Due to this special property, we refer to such an even cycle a
\textit{half-negatived even cycle}.

Hassani Monfare and Mallik [\citenum{even-cycle}]
raised a question of whether one can construct a specific matrix
for a given graph
that can be decomposed into a product of $RR^\tp$
where $R$ is another matrix associated with that graph.
This specific matrix is required to have a non-zero determinant
for even cycles and zero determinant for odd cycles.
However, we have demonstrated that
Hermitian Laplacian matrices cannot satisfy
those requirements no matter the choice of $\omega$.
Future work may explore alternative matrix constructions
to circumvent this limitation.

\section{Applications}\label{sec:applications}
We introduce some notations that will be frequently used in this section.
\begin{equation}
    \mathcal{C}(n,k,g)\coloneqq\left\{(V',E')\subseteq G:
    \substack{|V'|=n,\,(V',E')\text{ forms a }k\text{-unicyclic graph, whose cycle is:} \\
    \text{a }k\text{-cycle, with }g\text{ negative edges, with }n-k\text{ pendant edges attached}}\right\}.
\end{equation}
\begin{equation}
    \mathcal{C}(n,k,\cdot)\coloneqq\left\{(V',E')\subseteq G:
    \substack{|V'|=n,\,(V',E')\text{ forms a }k\text{-unicyclic graph, whose cycle is:} \\
    \text{a }k\text{-cycle, with arbitrarily many negative edges, with }n-k\text{ pendant edges attached}}\right\}.
\end{equation}
\begin{equation}
    \mathcal{F}(n)\coloneqq\left\{(V',E')\subseteq G:(V',E')\text{ is a rootless forest of }n\text{ vertices}\right\}.
\end{equation}

\subsection{Hermitian Laplacian Principal Minor Decomposition Using Fifth Roots of Unity}\label{ss:p5}
\subsubsection{Principal Minor Decomposition Under Parameter $p=5$}
The following application shows that
some specific choices of $\omega$
can effectively simplify the determinant expression,
enabling targeted analysis of certain types of substructures.

The choice of $p=5$
provides a special algebraic structure
for decomposing and interpreting
the Hermitian Laplacian determinants.
When $\omega=\omega_5\coloneqq\ee^{\im 2\pi/5}$,
the distinct non-zero determinant contributions
from a $k$-cycle with $g$ negative edges, i.e.
$\ell_{\omega_5}(k,g)=4\sin^2[(k-2g)\pi/5]$,
take on specific values within $\Q(\sqrt{5})$.

In fact, a unicyclic graph can only contribute
these three different real numbers:
\begin{equation}
    \alpha=4\sin^2\left(\frac{2\pi}{5}\right)
    =4\sin^2\left(\frac{3\pi}{5}\right)
    =\frac{5+\sqrt{5}}{2}\approx 3.618,
\end{equation}
\begin{equation}
    \beta=4\sin^2\left(\frac{\pi}{5}\right)
    =4\sin^2\left(\frac{4\pi}{5}\right)
    =\frac{5-\sqrt{5}}{2}\approx 1.382,
\end{equation}
and $0$, to the Hermitian Laplacian determinant.
These two numbers are in (and are linearly independent over)
the field generated by $\sqrt{5}$ over $\Q$,
denoted by $\Q(\sqrt{5})$:
\begin{equation}
    \alpha,\beta\in\Q(\sqrt{5})=\{a+b\sqrt{5}:a,b\in\Q\}.
\end{equation}

Observe that
$\alpha,\beta$ are strongly related to the
\textit{golden ratio} $\phi=(1+\sqrt{5})/2$, in fact,
\begin{equation}
\begin{aligned}
    \alpha&=\sqrt{5}\phi=\phi+2, \\
    \beta&=\sqrt{5}\phi^{-1}=\sqrt{5}(\phi-1), \\
    \alpha+\beta&=\alpha\beta=5, \\
    \alpha-\beta&=\sqrt{5}, \\
    \dfrac{\alpha}{\beta}&=\phi^2=\phi+1.
\end{aligned}
\end{equation}

By some calculation, we have
\begin{equation}
    \ell(k,g)=\begin{cases}
        \alpha=4\sin^2\left(\frac{2\pi}{5}\right), &k-2g\equiv 2,3\pmod 5, \\
        \beta=4\sin^2\left(\frac{\pi}{5}\right), &k-2g\equiv 1,4\pmod 5, \\
        0, &k-2g\equiv 0\pmod 5.
    \end{cases}
\end{equation}

\begin{definition}
    If a unicyclic graph $C\in\mathcal{C}(n,k,g)$ satisfies $\ell_{\omega_5}(k,g)=0$,
    which happens when and only when $5\mid k-2g$,
    we call $(k,g)$ a \textit{vanishing pair}
    and $C$ a \textit{vanishing unicyclic graph}.

    If $k-2g\equiv 2,3\pmod 5$, we call $(k,g)$ an \textit{alpha pair}
    and $C$ an \textit{alpha unicyclic graph};
    if $k-2g\equiv 1,4\pmod 5$, we call $(k,g)$ a \textit{beta pair}
    and $C$ a \textit{beta unicyclic graph}.
\end{definition}

The following theorem demonstrates
how the number of these unicyclic graphs
is reflected in the determinant expansion.

\begin{theorem}
Suppose an oriented graph $G(V,E)$
whose Hermitian Laplacian matrix is $L_\omega$ where $\omega=\omega_5$,
and $V'\subseteq V$, and $\{E'_i\}\subseteq 2^E$ are
its edge subsets such that every $(V',E'_i)$
is an all-regular substructure that contains
$n_\alpha^{(i)}$ alpha unicyclic graphs,
$n_\beta^{(i)}$ beta unicyclic graphs,
$n_0^{(i)}$ vanishing graphs
and arbitrarily many rootless trees.
Then
\begin{equation}
    \det L_\omega[V']=\sum_i \delta_0(n_0^{(i)})
    (\sqrt{5})^{n_\alpha^{(i)}+n_\beta^{(i)}}
    \phi^{n_\alpha^{(i)}-n_\beta^{(i)}}
\end{equation}
where $i$ traverses each all-regular substructure $(V',E'_i)$
formed by $V'$.

Here, \begin{equation}
    \delta_0(x)=\begin{cases}
        1, &x=0, \\
        0, &x\ne 0
    \end{cases}
\end{equation}
is a Kronecker delta.
\end{theorem}

\begin{proof}
By Hermitian Laplacian Matrix Decomposition Theorem we have
\begin{equation}
\begin{split}
    \det L_\omega[V']&=\sum_i \det L_\omega((V',E'_i)) \\
    &=\sum_i \delta_0(n_0^{(i)})
    \alpha^{n_\alpha^{(i)}}
    \beta^{n_\beta^{(i)}} \\
    &=\sum_i \delta_0(n_0^{(i)})
    (\sqrt{5}\phi)^{n_\alpha^{(i)}}
    (\sqrt{5}\phi^{-1})^{n_\beta^{(i)}} \\
    &=\sum_i \delta_0(n_0^{(i)})
    (\sqrt{5})^{n_\alpha^{(i)}+n_\beta^{(i)}}
    \phi^{n_\alpha^{(i)}-n_\beta^{(i)}},
\end{split}
\end{equation}
completing the proof.
\end{proof}

Additionally, if we consider the Hermitian Laplacian determinant
of only one single substructure $G'(V',E')$
that contains no vanishing unicyclic graphs,
$n_\alpha$ alpha unicyclic graphs,
and $n_\beta$ beta unicyclic graphs,
then the determinant equals to
\begin{equation}
    \alpha^{n_\alpha}\beta^{n_\beta}
    =(\sqrt{5})^{n_\alpha+n_\beta}\phi^{n_\alpha-n_\beta}
    \eqqcolon s(n_\alpha,n_\beta).
\end{equation}

Crucially, it can be proven that:
\begin{theorem}\label{thm:bijection}
    The value of
    \begin{equation}
        s(x,y)\coloneqq\alpha^x\beta^y=(\sqrt{5}\phi)^x(\sqrt{5}\phi^{-1})^y
        =(\sqrt{5})^{x+y}\phi^{x-y}
    \end{equation}
    corresponds bijectively to the
    non-negative integer pair $(x,y)$.
\end{theorem}
\begin{proof}
    Write $x+y=A,x-y=B$. Then $s(x,y)=(\sqrt{5})^A\phi^B$.
    For the sake of contradiction,
    assume that there exists another non-negative integer pair $(x',y')$
    such that $s(x,y)=s(x',y')$,
    and write $x'+y'=A',x'-y'=B'$.
    
    Write $A-A'=C,B-B'=D$. Then
    \begin{equation}\label{eq:bij_proof_1}
    \begin{aligned}
        &s(x,y)=(\sqrt{5})^A\phi^B=(\sqrt{5})^{A'}\phi^{B'}=s(x',y') \\
        \Rightarrow& (\sqrt{5})^{A-A'}\phi^{B-B'}=(\sqrt{5})^C\phi^D=1.
    \end{aligned}
    \end{equation}

    Consider field norm $\mathcal{N}:\Q(\sqrt{5})\to\Q,\,a+b\sqrt{5}\mapsto a^2-5b^2$.
    Then $\mathcal{N}$ satisfies that
    \begin{itemize}
        \item $\mathcal{N}(uv)=\mathcal{N}(u)\mathcal{N}(v)$ for $u,v\in\Q(\sqrt{5})$;
        \item $\mathcal{N}(u^n)=(\mathcal{N}(u))^n$ for $u\in\Q(\sqrt{5}),n\in\Z$.
    \end{itemize}

    Taking the field norms for both sides in Eq. (\ref{eq:bij_proof_1})
    we have
    \begin{equation}
    \begin{aligned}
        &\mathcal{N}\left((\sqrt{5})^C\phi^D\right)=\mathcal{N}(1) \\
        \Rightarrow& \left(\mathcal{N}(\sqrt{5})\right)^C \left(\mathcal{N}(\phi)\right)^D=1 \\
        \Rightarrow& (-5)^C (-1)^D=1 \\
        \Rightarrow& C=0,
    \end{aligned}
    \end{equation}
    therefore
    \begin{equation}
        1=(\sqrt{5})^0\phi^D=\phi^D\Rightarrow D=0.
    \end{equation}

    Since $C=D=0$, we have $A=A',B=B'$.
    Observe that $x,y$ and $x',y'$ are uniquely determined by
    $A,B$ and $A',B'$, we have $x=x',y=y'$.
    This completes the proof.
\end{proof}

If one wants to utilize this theorem
to deduce $n_\alpha$ and $n_\beta$ from
the Hermitian Laplacian determinant
of a certain all-regular substructure,
they have two approaches:
\begin{enumerate}[label=(\arabic*)]
    \item use symbolic computation
    to calculate this determinant
    and obtain the exact value $a+b\sqrt{5}$,
    then perform operations on $a$ and $b$
    to get $n_\alpha$ and $n_\beta$; or
    \item use numerical computation
    to calculate the determinant
    and obtain a numerical value,
    then look up this value in a table.
\end{enumerate}
The disadvantage of (1) is
that symbolic computation is typically much slower
than numerical computation;
the disadvantage of (2) is
that numerical values inevitably carry errors,
and the table lookup process requires additional time.

However, in some cases, we care about $n_\alpha\pm n_\beta$
more than $n_\alpha$ and $n_\beta$ themselves.
In fact, we can solve for $n_\alpha$ and $n_\beta$
from $n_\alpha+n_\beta$ and $n_\alpha-n_\beta$.
In the next section, we propose a theorem that
gives an exact formula for $n_\alpha\pm n_\beta$
at the cost of one additional computation
of the Hermitian Laplacian determinant
under parameter $p=5/2$.

\subsubsection{A Unicyclic Graph Enumeration Method Using Parameters $p=5$ and $p=5/2$}
\label{sec:5-and-5/2}
Let $\zeta=\ee^{\im 2\pi/5}$, then
\begin{equation}
\begin{aligned}
    \Q(\sqrt{5})&=\Q(\alpha)=\Q(\beta) \\
    &=\Q\left(4\sin^2\left(\dfrac{\pi}{5}\right)\right) \\
    &=\Q\left(2-2\cos\left(\dfrac{2\pi}{5}\right)\right) \\
    &=\Q\left(\cos\left(\dfrac{2\pi}{5}\right)\right)
\end{aligned}
\end{equation}
is the maximal real subfield of $\Q(\zeta)$.
In fact,
\begin{equation}
\begin{aligned}
    \cos\left(\dfrac{2\pi}{5}\right)&=\dfrac{1}{2}(\zeta+\zeta^{-1}), \\
    \sqrt{5}&=\zeta+\zeta^4-\zeta^2-\zeta^3.
\end{aligned}
\end{equation}

We want to find an automorphism
$\sigma\in\gal(\Q(\zeta)/\Q)$
whose restriction to the subfield $\Q(\sqrt{5})$
is the non-trivial automorphism of $\gal(\Q(\sqrt{5})/\Q)$,
which exists by Galois theory.
That is, we require $\sigma(\sqrt{5})=-\sqrt{5}$.
Observe that $\sigma_2:\zeta\mapsto\zeta^2$
is an automorphism and satisfies:
\begin{equation}
\begin{aligned}
    \sigma_2(\sqrt{5})&=\sigma_2\left(\zeta+\zeta^4-\zeta^2-\zeta^3\right) \\
    &=\zeta^2+\zeta^8-\zeta^4-\zeta^6 \\
    &=\zeta^2+\zeta^3-\zeta^4-\zeta \\
    &=-\sqrt{5},
\end{aligned}
\end{equation}
thus $\sigma_2$ is the desired automorphism.

By the discussion above, we have the following proposition.
\begin{proposition}
    Let $\omega_5=\ee^{\im 2\pi/5}$,
    $\omega_5^2=\ee^{\im 4\pi/5}=\ee^{\im 2\pi/(5/2)}=\omega_{5/2}$.
    Then $\det L_{\omega_5}$ and $\det L_{\omega_5^2}$
    are Galois conjugates, i.e.,
    \begin{equation}
        \det L_{\omega_5}=a+b\sqrt{5}\quad\Leftrightarrow\quad
        \det L_{\omega_5^2}=a-b\sqrt{5}.
    \end{equation}
    Furthermore,
    \begin{equation}
        \det L_{\omega_5}=\alpha^{n_\alpha}\beta^{n_\beta}\quad\Leftrightarrow\quad
        \det L_{\omega_5^2}=\beta^{n_\alpha}\alpha^{n_\beta}.
    \end{equation}
\end{proposition}
\begin{proof}
    Suppose $\det L_{\omega_5}=\alpha^{n_\alpha}\beta^{n_\beta}$. Then
    \begin{equation}
    \begin{aligned}
        \det L_{\omega_5^2}&=\det(\sigma_2(L_{\omega_5})) \\
        &=\sigma_2\left(\det L_{\omega_5}\right) \\
        &=\sigma_2(\alpha^{n_\alpha}\beta^{n_\beta}) \\
        &=(\sigma_2(\alpha))^{n_\alpha}(\sigma_2(\beta))^{n_\beta} \\
        &=\beta^{n_\alpha}\alpha^{n_\beta}
    \end{aligned}
    \end{equation}
    and vice versa.

    Here, for some matrix $A$ whose elements are in $\Q(\zeta)$,
    $\sigma_2(A)$ denotes the new matrix formed by
    applying $\sigma_2$ to every element of $A$.
\end{proof}

By this proposition we can obtain the exact expression
for $n_\alpha\pm n_\beta$, stated as the theorem below.
\begin{theorem}
    Let $\omega_5=\ee^{\im 2\pi/5}$,
    $\omega_5^2=\ee^{\im 4\pi/5}=\ee^{\im 2\pi/(5/2)}=\omega_{5/2}$.
    Suppose $G'(V',E')$ is an all-regular substructure
    that contains
    $n_\alpha$ alpha unicyclic graphs,
    $n_\beta$ beta unicyclic graphs
    and $n_0=0$ vanishing unicyclic graphs.
    Denote $L_{\omega_5}$ and $L_{\omega_5^2}$
    to be its two Hermitian Laplacian matrices,
    parameterized by $\omega_5$ and $\omega_5^2$ respectively.
    Then we have
    \begin{equation}
        \begin{cases}
            n_\alpha+n_\beta=\log_5[(\det L_{\omega_5})(\det L_{\omega_5^2})], \\
            n_\alpha-n_\beta=\log_{\phi^2}\left(\dfrac{\det L_{\omega_5}}{\det L_{\omega_5^2}}\right).
        \end{cases}
    \end{equation}
\end{theorem}
\begin{proof}
    By direct calculation, we have
    \begin{equation}
    \begin{aligned}
        (\det L_{\omega_5})(\det L_{\omega_5^2})
        &=\alpha^{n_\alpha}\beta^{n_\beta}\cdot\beta^{n_\alpha}\alpha^{n_\beta} \\
        &=\alpha^{n_\alpha}\beta^{n_\alpha}\cdot\alpha^{n_\beta}\beta^{n_\beta} \\
        &=(\alpha\beta)^{n_\alpha+n_\beta}=5^{n_\alpha+n_\beta}.
    \end{aligned}
    \end{equation}
    and
    \begin{equation}
    \begin{aligned}
        \dfrac{\det L_{\omega_5}}{\det L_{\omega_5^2}}
        &=\dfrac{\alpha^{n_\alpha}\beta^{n_\beta}}{\beta^{n_\alpha}\alpha^{n_\beta}} \\
        &=\left(\dfrac{\alpha}{\beta}\right)^{n_\alpha}
        \left(\dfrac{\beta}{\alpha}\right)^{n_\beta} \\
        &=\left(\dfrac{\alpha}{\beta}\right)^{n_\alpha}
        \left(\dfrac{\alpha}{\beta}\right)^{-n_\beta} \\
        &=\left(\dfrac{\alpha}{\beta}\right)^{n_\alpha-n_\beta}
        =(\phi^2)^{n_\alpha-n_\beta},
    \end{aligned}
    \end{equation}
    completing the proof.
\end{proof}

\subsection{A Unicyclic Graph Enumeration Method by
Parameterizing Hermitian Laplacian Matrices
on Galois Conjugates}
This section serves as a generalization of Section \ref{sec:5-and-5/2}.
In that section, we demonstrated that when $p=5$,
we can take the Galois conjugates of $\omega_5$,
namely $\omega_{5/2}$, as the new parameter,
then multiplying and dividing the original determinant
by the new determinant
directly yields a power of $5$ and $\phi^2$.
In this section, we will prove that
the multiplicative case is not an exception:
for all odd prime number $p$,
the product of the Hermitian Laplacian determinants
of parameters $\omega_p$ and all its Galois conjugates
must be an integer power of $p$.

First, we present a theorem and a lemma
to lay the groundwork for deriving the key conclusion of this section.

\begin{theorem}[Sine Product Formula]\label{thm:sine-prod}
    Suppose $n\in\Z$ and $n\geqslant 2$. Then
    \begin{equation}
        \prod_{k=1}^{n-1}\sin\left(\dfrac{k\pi}{n}\right)=\dfrac{n}{2^{n-1}}.
    \end{equation}
\end{theorem}
\begin{proof}
    Consider the polynomial $z^n-1=0$,
    whose roots are $\omega_k=\ee^{\im 2\pi k/n}\,(k=0,1,\ldots,n-1)$.
    We have
    \begin{equation}
        z^n-1=\prod_{k=0}^{n-1}(z-\omega_k).
    \end{equation}
    Removing the root $z=1$ (corresponding to $k=0$), we have
    \begin{equation}
        \dfrac{z^n-1}{z-1}=\prod_{k=1}^{n-1}(z-\omega_k).
    \end{equation}
    This implies that
    \begin{equation}
        \lim_{z\to 1}\dfrac{z^n-1}{z-1}=\lim_{z\to 1}\prod_{k=1}^{n-1}(z-\omega_k)
    \end{equation}
    \iffalse
    where
    \begin{equation}
    \begin{aligned}
        \text{LHS}&=\left.\dfrac{\mathrm{d}}{\mathrm{d}z}z^n\right|_{z=1}
        =\left.nz^{n-1}\right|_{z=1}&=n, \\
        \text{RHS}&=\prod_{k=1}^{n-1}(1-\omega_k),
    \end{aligned}
    \end{equation}
    \fi
    therefore
    \begin{equation}
        \prod_{k=1}^{n-1}(1-\omega_k)=n.
    \end{equation}
    Take the modulus of both sides:
    \begin{equation}
    \begin{aligned}
        \text{LHS}&=\prod_{k=1}^{n-1}|1-\omega_k|
        =\prod_{k=1}^{n-1}|1-\ee^{\im 2\pi k/n}| \\
        &=\prod_{k=1}^{n-1}\left|\ee^{\im\pi k/n}(\ee^{-\im\pi k/n}-\ee^{\im\pi k/n})\right|=\prod_{k=1}^{n-1}\left|\ee^{\im\pi k/n}\left(-2\im\cdot\sin\left(\dfrac{k\pi}{n}\right)\right)\right| \\
        &=\prod_{k=1}^{n-1}|\ee^{\im\pi k/n}|\cdot|-2\im|\cdot\left|\sin\left(\dfrac{k\pi}{n}\right)\right|=\prod_{k=1}^{n-1}2\left|\sin\left(\dfrac{k\pi}{n}\right)\right| \\
        &=2^{n-1}\prod_{k=1}^{n-1}\sin\left(\dfrac{k\pi}{n}\right)=\text{RHS}=n,
    \end{aligned}
    \end{equation}
    which is what we want.
\end{proof}

\begin{corollary}
    Suppose $p\geqslant 3$ is a positive odd number, then
    \begin{equation}
        \prod_{n=1}^{(p-1)/2}\left[4\sin^2\left(\dfrac{n\pi}{p}\right)\right]=p.
    \end{equation}
\end{corollary}
\begin{proof}
    By Theorem \ref{thm:sine-prod} and notice that
    \begin{equation}
        \sin\left(\dfrac{n\pi}{p}\right)=\sin\left[\dfrac{(p-n)\pi}{p}\right],
    \end{equation}
    we have
    \begin{equation}
        \prod_{n=1}^{(p-1)/2}\sin\left(\dfrac{n\pi}{p}\right)
        =\left[\prod_{n=1}^{p-1}\sin\left(\dfrac{n\pi}{p}\right)\right]^{1/2}
        =\left(\dfrac{p}{2^{p-1}}\right)^{1/2}.
    \end{equation}
    Therefore
    \begin{equation}
    \begin{aligned}
        \prod_{n=1}^{(p-1)/2}\left[4\sin^2\left(\dfrac{n\pi}{p}\right)\right]
        =4^{(p-1)/2}\left[\prod_{n=1}^{(p-1)/2}\sin\left(\dfrac{n\pi}{p}\right)\right]^2
        =2^{p-1}\cdot\left[\left(\dfrac{p}{2^{p-1}}\right)^{1/2}\right]^2
        =p,
    \end{aligned}
    \end{equation}
    completing the proof.
\end{proof}

\begin{lemma}\label{lem:latin-square}
    Let $p$ be an odd prime number
    and let $S=\{1,2,\ldots,(p-1)/2\}$.
    We define a binary operation $\ast$ on $S$ as follows.
    
    For any $a, x \in S$:
    \begin{equation}
        a\ast x=\begin{cases}
            ax \bmod p, & \text{if } 1\leqslant(ax\bmod p)\leqslant\frac{p-1}{2}, \\
            p-(ax\bmod p), & \text{if } \frac{p-1}{2}<(ax\bmod p)\leqslant p-1.
        \end{cases}
    \end{equation}
    
    Then the Cayley table for $(S,\ast)$ forms a Latin Square.
    This means that for any fixed element $a\in S$,
    the functions $f_a(x)=a\ast x$ (representing the values on the rows)
    and $g_a(x)=x\ast a$ (representing the values on the columns)
    are both permutations of $S$.
\end{lemma}

\begin{proof}
    We show that each element of $S$ appears exactly once
    in each row and each column of the Cayley table.
    It is sufficient to prove that
    for every $a\in S$, both $f_a$ and $g_a$ are injective.
    Because the operation $\ast$ is commutative
    since the product $ax$ is commutative,
    the Cayley table is symmetric.
    Thus we only need to prove the injectivity for $f_a$;
    the injectivity for $g_a$ will follow automatically.

    Let $a\in S$ be a fixed element.
    For the sake of contradiction, assume that
    for some $x_1,x_2\in S,x_1\ne x_2$,
    we have $f_a(x_1)=f_a(x_2)$.

    By the definition of the operation $\ast$,
    the condition $a\ast x_1=a \ast x_2$ implies that
    $ax_1\bmod p$ and $ax_2\bmod p$ are either equal or sum to $p$.
    This can be expressed as:
    \begin{equation}
        ax_1 \equiv\pm ax_2 \pmod{p}.
    \end{equation}
    This gives rise to two cases.

    \paragraph{Case 1: $ax_1\equiv ax_2\pmod{p}$.}
    Since $p$ is a prime number and $a<p$,
    we have $\gcd(a,p)=1$.
    Thus $a$ has a multiplicative inverse modulo $p$
    and we can cancel $a$ from the congruence:
    \begin{equation}
        x_1\equiv x_2\pmod{p}.
    \end{equation}
    Therefore $(x_1-x_2)\bmod p=0$.
    However, since $1\leqslant x_1, x_2\leqslant\frac{p-1}{2}$,
    the difference $|x_1-x_2|$ must be less than $p$.
    Therefore $x_1-x_2=0$, which is a contradiction.

    \paragraph{Case 2: $ax_1\equiv -ax_2\pmod{p}$.}
    Similarly, we can cancel $a$ since $\gcd(a,p) = 1$:
    \begin{equation}
        x_1\equiv -x_2\pmod{p}.
    \end{equation}
    Therefore $(x_1+x_2)\bmod p=0$,
    which means that $x_1+x_2$ is a multiple of $p$.
    However, the sum is bounded:
    $2\leqslant x_1+x_2\leqslant 2\cdot\frac{p-1}{2}=p-1$.
    So this is also a contradiction.
    
    This completes the proof.
\end{proof}

Similar to the Section \ref{sec:5-and-5/2},
as a generalization of the vanishing unicyclic graph,
we introduce the notion of $p$-vanishing unicyclic graph.

\begin{definition}
    For a given parameter $\omega=\omega_p=\ee^{\im 2\pi/p}$,
    the Hermitian Laplacian determinant contribution from a cycle is
    $\ell_{\omega_p}(k,g)=4\sin^2[(k-2g)\pi/p]=2-2\cos[2(k-2g)\pi/p]$,
    where $k$ is the cycle length and $g$ is the number of negative edges.
    This value is non-zero if and only if $p\nmid (k-2g)$.
    We refer to such components as \textit{non-$p$-vanishing},
    and as \textit{$p$-vanishing} otherwise.
\end{definition}

With the preparatory work completed,
we now present the main result of this section.

\begin{theorem}
    Let $p$ be an odd prime number and denote
    $\omega_{p/q}=\ee^{\im 2\pi/(p/q)}\,(q\in\{0,1,\ldots,p-1\})$
    to be the primitive $p$th roots of unity.
    Suppose $G'(V',E')$ is an all-regular substructure.
    Denote $\{L_{\omega_{p/q}}\}_{q=1}^{(p-1)/2}$
    to be a set of its Hermitian Laplacian matrices,
    each of which is parameterized by $\omega_{p/q}$.
    Suppose $G'$ contains no $p$-vanishing unicyclic graphs
    and denote $n_\ast$ to be the total number of
    non-$p$-vanishing unicyclic graphs.
    Then we have
    \begin{equation}
        n_\ast=\log_p\left(\prod_{q=1}^{(p-1)/2}\det L_{\omega_{p/q}}\right).
    \end{equation}
\end{theorem}

\begin{proof}
    Let the distinct non-zero cycle determinant values
    for the parameter $\omega_p$ be
    $\ell_1,\ell_2,\ldots,\ell_m$, where $m=(p-1)/2$
    and $\ell_j=4\sin^2(j\pi/p)$ for $j\in\{1,\ldots,m\}$.
    Assume $n_j$ is the number of unicyclic graphs $C\in\mathcal{C}(\cdot,k,g)$
    such that $k-2g\equiv j\,(\tmod\,p)$.
    Then such unicyclic graphs have
    Hermitian Laplacian determinant $\ell_j$.
    Therefore
    \begin{equation}
        \det L_{\omega_p}(G')=\prod_{j=1}^{m}\ell_j^{n_j}.
    \end{equation}

    The set of values $\{\cos(2j\pi/p)\}_{j=1}^m$
    generates the maximal real subfield of the $p$th cyclotomic field $\Q(\omega_p)$,
    which we denote by $K_p=\Q(\cos(2\pi/p))$.
    By Galois theory, $K_p/\Q$ is a Galois extension of degree $m=(p-1)/2$
    and is a Galois subextension of $\Q(\omega_p)/\Q$;
    and the automorphisms of the Galois group $\gal(K_p/\mathbb{Q})$
    are the automorphisms
    $\sigma_p\in\gal(\Q(\omega_p)/\Q):\omega_p\mapsto\omega_p^q$
    restricted to $K_p$.
    In fact, it holds that
    \begin{equation}
        \begin{aligned}
        \sigma_q\left(\cos\left(\dfrac{2j\pi}{p}\right)\right)
        &=\sigma_q\left(\dfrac{\ee^{2j\pi/p}+\ee^{-2j\pi/p}}{2}\right) \\
        &=\dfrac{\sigma_q(\ee^{2j\pi/p})+\sigma_q(\ee^{-2j\pi/p})}{2} \\
        &=\dfrac{\ee^{2jq\pi/p}+\ee^{-2jq\pi/p}}{2}
        =\cos\left(\dfrac{2jq\pi}{p}\right).
        \end{aligned}
    \end{equation}
    Specifically,
    it shows that $\sigma_q$ acts on the elements $\ell_j\in K_p$ by
    \begin{equation}
        \sigma_q(\ell_j)=\sigma_q\left(4\sin^2\left(\dfrac{j\pi}{p}\right)\right)
        =4\sin^2\left(\dfrac{jq\pi}{p}\right)=\ell_{[jq]_p}.
    \end{equation}

    Here, the notation $[n]_p$ is defined as\footnote{
    For example, we have $[2]_7=[5]_7=[9]_7=[12]_7$, and therefore
    $4\sin^2(2\pi/7)=4\sin^2(5\pi/7)=4\sin^2(9\pi/7)=4\sin^2(12\pi/7).$
    }
    \begin{equation}
        [n]_p=\begin{cases}
            n\bmod p, & \text{if } 1\leqslant(n\bmod p)\leqslant\frac{p-1}{2}, \\
            p-(n\bmod p), & \text{if } \frac{p-1}{2}<(n\bmod p)\leqslant p-1.
        \end{cases}
    \end{equation}

    %Back to the topic.
    We now consider the determinant of the matrix $L_{\omega_{p/q}}$.
    Under parameter $\omega_{p/q}$,
    the Hermitian Laplacian determinant
    of a unicyclic graph $C\in\mathcal{C}(\cdot,k,g)$
    such that $k-2g\equiv j\,(\tmod\,p)$ is
    \begin{equation}
        4\sin^2\left(\dfrac{j\pi}{p/q}\right)
        =4\sin^2\left(\dfrac{jq\pi}{p}\right)
        =\sigma_q\left(4\sin^2\left(\dfrac{j\pi}{p}\right)\right),
    \end{equation}
    which is exactly applying the automorphism $\sigma_q$
    to the determinant under parameter $\omega_p$.
    Therefore we have
    \begin{equation}
        \det L_{\omega_{p/q}}(G')
        =\sigma_q\left(\det L_{\omega_p}(G')\right)
        =\sigma_q\left(\prod_{j=1}^{m}\ell_j^{n_j}\right)
        =\prod_{j=1}^{m}(\sigma_q(\ell_j))^{n_j}.
    \end{equation}

    Let's compute the product of $\det L_{\omega_{p/q}}(G')$
    over all $q\in\{1,\ldots,m\}$:
    \begin{equation}
    \begin{aligned}
        \prod_{q=1}^{m}\det L_{\omega_{p/q}}&=\prod_{q=1}^{m}\left(\prod_{j=1}^{m}(\sigma_q(\ell_j))^{n_j}\right) \\
        &=\prod_{j=1}^{m}\left(\prod_{q=1}^{m}(\sigma_q(\ell_j))^{n_j}\right) \\
        &=\prod_{j=1}^{m}\left(\prod_{q=1}^{m}\sigma_q(\ell_j)\right)^{n_j}.
    \end{aligned}
    \end{equation}

    By Lemma \ref{lem:latin-square},
    it holds that $(\sigma_q(\ell_1),\ldots,\sigma_q(\ell_m))$
    is a permutation of $(\ell_1,\ldots,\ell_m)$.
    Therefore, for any fixed $j$, the set $\{\sigma_q(\ell_j)\}_{q=1}^m$
    is the set of all Galois conjugates of $\ell_j$,
    which is the set $\{\ell_j\}_{j=1}^m$ itself.
    Therefore
    \begin{equation}
        \prod_{q=1}^{m}\sigma_q(\ell_j)=\prod_{j=1}^{m}\ell_j.
    \end{equation}
    Moreover, by Theorem \ref{thm:sine-prod},
    this product is constant for all $j$:
    \begin{equation}
        \prod_{j=1}^{m}\ell_j=p.
    \end{equation}

    Substituting these back into our expression, we have
    \begin{equation}
        \prod_{q=1}^{m}\det L_{\omega_{p/q}}=\prod_{j=1}^{m} p^{n_j}=p^{\sum_{j=1}^{m} n_j}=p^{n_\ast}.
    \end{equation}
    The theorem follows
    directly by taking the logarithm base $p$ of this result.
\end{proof}

\begin{remark}
    The reader may wonder why the condition
    of $p$ being an odd prime is essential for this theorem.
    We discuss the consequences of relaxing this condition
    to any odd integer, using $p=9$ as a counterexample.

    When $p=9$, the structure of the Galois group $\gal(\Q(\omega_9)/\Q)$
    becomes different:
    its elements are $\sigma_n:\omega_9\mapsto\omega_9^n$
    where $n$ is coprime to $p$.
    For example, the map $\sigma_3:\omega_9\mapsto\omega_9^3$
    is not an automorphism on $\Q(\omega_9)$ since $\gcd(3,9)\ne 1$,
    and $\omega_9^3=\omega_3$ is not a primitive $9$th root of unity.
    This results in any further mapping $\sigma_n(\omega_9^3)$
    only yields element in $\Q(\omega_3)\subset\Q(\omega_9)$
    rather than the whole $\Q(\omega_9)$.
\end{remark}

\subsection{Some Relationships Among Hermitian Laplacian Principal Minors Under Parameters $p=2,3,4,5,6$}
We take $\omega_2=-1,\,\omega_3=\ee^{\im 2\pi/3}=(-1+\sqrt{3}\im)/2,\,\omega_4=\im,\,\omega_5=\ee^{\im 2\pi/5}$
and $\omega_6=\ee^{\im 2\pi/6}=(1+\sqrt{3}\im)/2$ in this section.

In this section, we decompose
the Hermitian Laplacian principal minors
under different parameters $\omega_i\,(i\in\{2,3,4,5,6\})$
and derive some relationships among these values
as an application of this paper.
Furthermore, these relationships may potentially be used
to detect the presence of triangles and $4$-cycles
within given vertex subsets.

We introduce the following notations:
\begin{equation}
    \mathcal{TU}(n,k,g)\coloneqq\left\{(V',E')\subseteq G:
    \substack{(V',E')\text{ forms a TU-substructure,}\\\text{which is composed of one rootless forest }\mathcal{F}(\cdot)\text{ and one }k\text{-unicyclic graph }\mathcal{C}(n,k,g)}\right\}.
\end{equation}
\begin{equation}
    \# S\coloneqq\text{size of finite set }S=|S|.
\end{equation}

\subsubsection{On Triangles}
Let's say we want to decompose the Hermitian Laplacian principal minor
corresponding to the vertex subset $V'=\{v_1,v_2,v_3\}$,
which can possibly form a triangle (a $3$-cycle).

\begin{equation}
    \displaystyle\det L_\omega[V']=\sum_{\substack{E'\subseteq E\\|E'|=|V'|}}\det L_\omega(G'(V',E')).
\end{equation}

We take $\omega=\omega_4=\im$ first.
Let $E'$ be an arbitrary edge subset such that $|E'|=3$.
If $(V',E')$ forms a triangle $\mathcal{C}(3,3,\cdot)$,
then by Proposition \ref{prp:cycle_lap},
it contributes
\begin{equation}
    \begin{aligned}
        \ell_{\omega_4}(3,g)&=4\sin^2\left[\frac{(3-2g)\pi}{4}\right]\quad(g\in\{0,1\}) \\
        &=4\sin^2\left(\frac{\pi}{4}\right) \\
        &=2
    \end{aligned}
\end{equation}
to $\det L_\omega[V']$, therefore we have
\begin{equation}
    \begin{aligned}
        \displaystyle\det L_{\omega_4}[V']
        &=\quad\#\{\mathcal{F}(3)\text{ spanned by }V'\} \\
        &\quad +2\cdot\#\{\mathcal{C}(3,3,\cdot)\text{ spanned by }V'\}, \\
    \end{aligned}
\end{equation}
i.e. $\det L_{\omega_4}[V']$ equals to
\begin{equation}\label{eq:triangle-omega4}
    \text{number of rootless trees spanned by }V'+2\times\text{number of triangles spanned by }V'.
\end{equation}

We take $\omega=\omega_2=-1$ next.
By generalized Matrix-Tree Theorem for mixed graphs
by Bapat [\citenum{bapat}],
or by some similar analysis using
the parameterized Hermitian Laplacian framework,
we have
\begin{equation}
    \det L_{-1}[V']=\sum_{\substack{E'\subseteq E\\(V',E')\text{ all-regular}}}4^{\#\mathcal{C}(\cdot,\text{odd},\cdot)},
\end{equation}
where $\#\mathcal{C}(\cdot,\text{odd},\cdot)$ denotes
the number of odd unicyclic graph components in $(V',E')$.

Additionally, similar to the approach employed by
Hassani Monfared and Mallik in [\citenum{even-cycle}]
to derive an analog of the matrix tree theorem for
the signless Laplacian matrix,
if we adopt the concept of \textit{TU-graph}
introduced by Cvetkovi\'c et al. in [\citenum{tu-graph}]
and suppose every all-regular substructure spanned by
$V'$ is a TU-substructure,
i.e. its connected components consist of only
rootless trees and at most $1$ unicyclic graph,
the above expression can be substantially simplified into:

\begin{equation}\label{eq:triangle-omega2-tu}
    \begin{aligned}   
    \det L_{-1}[V']
    &=\quad\#\{\mathcal{F}(\cdot)\text{ spanned by }V'\} \\
    &\quad+4\cdot\#\{\mathcal{TU}(\cdot,\text{odd},\cdot)\text{ spanned by }V'\}.
    \end{aligned}
\end{equation}

Since an all-regular substructure spanned by
a vertex subset with no more than $5$ vertices
can only be a TU-substructure,
we have $\mathcal{TU}(3,k,g)=\mathcal{C}(3,k,g)$.
Therefore

\begin{equation}\label{eq:triangle-omega2}
    \begin{aligned}   
    \det L_{-1}[V']
    &=\quad\#\{\mathcal{F}(3)\text{ spanned by }V'\} \\
    &\quad+4\cdot\#\{\mathcal{C}(3,3,\cdot)\text{ spanned by }V'\} \\
    \end{aligned}
\end{equation}

By comparing Eq. (\ref{eq:triangle-omega4})
and Eq. (\ref{eq:triangle-omega2}),
we have the following theorem.

\begin{theorem}
    Let $G(V,E)$ be a graph,
    $L_\omega$ be its Hermitian Laplacian matrix
    and $V'\subseteq V$ be its vertex subset
    such that $|V'|=3$,
    then we have
    \begin{equation}
        \begin{aligned}
            &\quad \det L_{\omega_2}[V']-\det L_{\omega_4}[V'] \\
            &=2\times\text{the number of triangles spanned by }V'.
        \end{aligned}
    \end{equation}
    and
    \begin{equation}
        \begin{aligned}
            &\quad 2\det L_{\omega_4}[V']-\det L_{\omega_2}[V'] \\
            &=\text{the number of rootless trees spanned by }V'.
        \end{aligned}
    \end{equation}
\end{theorem}

\subsubsection{On $4$-cycles}
The case of 3-cycles is relatively simple,
we now proceed to examine 4-cycles.
We demonstrate another application
of the parameterized Hermitian Laplacian framework:
the precise algebraic enumeration of distinct types
of $4$-vertex all-regular substructures
spanned by a given vertex subset
$V'=\{v_1,v_2,v_3,v_4\}$.

A $4$-vertex all-regular substructure can be categorized
into the following cases:

\begin{enumerate}[label=(\arabic*)]
    \item One $4$-cycle;
    \item One $3$-cycle and one $1$-rootless tree, which is an odd TU-substructure;
    \item One $3$-unicyclic graph of $4$ vertices, which is also an odd TU-substructure;
    \item One rootless forest, i.e. all components are rootless trees,
    amounting to a total of $p(4)=5$ distinct types
    (corresponding to the partitions $4$, $3+1$, $2+2$, $2+1+1$, and $1+1+1+1$).
    Here, $p(n)$ is the integer partition function.
\end{enumerate}

\begin{remark}
    The odd TU-substructures mentioned above, i.e.
    $\mathcal{TU}(3,3,g)\,(g=0,1)$,
    consist of exactly such two types:
    either (2) $\mathcal{C}(3,3,g)\sqcup\mathcal{F}(1)$
    or (3) $\mathcal{C}(4,3,g)$.
\end{remark}

Therefore
$\mathcal{C}(4,4,0)$, $\mathcal{C}(4,4,1)$,
$\mathcal{TU}(3,3,0)$, $\mathcal{TU}(3,3,1)$, and $\mathcal{F}(4)$
represent all the distinct types of all-regular substructures
which make distinct contributions to $\det L_\omega[V']$.
Our goal is to determine the numbers of these substructure
by choosing the parameter $\omega$ to be
some appropriate roots of unity.

\paragraph{$p=2$.}
By Eq. (\ref{eq:triangle-omega2-tu}) and
\begin{equation}
    k=4\quad\Rightarrow\quad
    \dfrac{k-2g}{2}\in\Z\,(\forall g\in\{0,1,2\})\quad\Rightarrow\quad
    \ell(4,g)\equiv 0
\end{equation}
we have
\begin{equation}
    \begin{aligned}
        \det L_{\omega_2}[V']&=\quad1\cdot\#\{\mathcal{F}(4)\text{ spanned by } V'\} \\
        &\quad +4\cdot\#\{\mathcal{TU}(\cdot,3,\cdot)\text{ spanned by } V'\}.
    \end{aligned}
\end{equation}

\paragraph{$p=3,4,6$.}
Similar to our discussions about triangles above,
we obtain that

\begin{equation}
\begin{aligned}
    \det L_{\omega_3}[V'] &= \quad 3\cdot\#\{\mathcal{C}(4,4,0)\text{ spanned by }V'\} \\
    &\quad +3\cdot\#\{\mathcal{C}(4,4,1)\text{ spanned by } V'\} \\
    &\quad +0\cdot\#\{\mathcal{C}(4,4,2)\text{ spanned by } V'\} \\
    &\quad +0\cdot\#\{\mathcal{C}(4,3,0)\text{ spanned by } V'\} \\
    &\quad +3\cdot\#\{\mathcal{C}(4,3,1)\text{ spanned by } V'\} \\
    &\quad +0\cdot\#\{\mathcal{C}(3,3,0)\sqcup\mathcal{F}(1)\text{ spanned by } V'\} \\
    &\quad +3\cdot\#\{\mathcal{C}(3,3,1)\sqcup\mathcal{F}(1)\text{ spanned by } V'\} \\
    &\quad +1\cdot\#\{\mathcal{F}(4)\text{ spanned by } V'\},
\end{aligned}
\end{equation}

\begin{equation}
\begin{aligned}
    \det L_{\omega_4}[V']&=\quad 0\cdot\#\{\mathcal{C}(4,4,0)\text{ spanned by }V'\} \\
    &\quad +4\cdot\#\{\mathcal{C}(4,4,1)\text{ spanned by }V'\} \\
    &\quad +0\cdot\#\{\mathcal{C}(4,4,2)\text{ spanned by }V'\} \\
    &\quad +2\cdot\#\{\mathcal{C}(4,3,0)\text{ spanned by }V'\} \\
    &\quad +2\cdot\#\{\mathcal{C}(4,3,1)\text{ spanned by }V'\} \\
    &\quad +2\cdot\#\{\mathcal{C}(3,3,0)\sqcup\mathcal{F}(1)\text{ spanned by }V'\} \\
    &\quad +2\cdot\#\{\mathcal{C}(3,3,1)\sqcup\mathcal{F}(1)\text{ spanned by }V'\} \\
    &\quad +1\cdot\#\{\mathcal{F}(4)\text{ spanned by }V'\},
\end{aligned}
\end{equation}

\begin{equation}
\begin{aligned}
    \det L_{\omega_6}[V']&=\quad 3\cdot\#\{\mathcal{C}(4,4,0)\text{ spanned by }V'\} \\
    &\quad +3\cdot\#\{\mathcal{C}(4,4,1)\text{ spanned by } V'\} \\
    &\quad +0\cdot\#\{\mathcal{C}(4,4,2)\text{ spanned by } V'\} \\
    &\quad +4\cdot\#\{\mathcal{C}(4,3,0)\text{ spanned by } V'\} \\
    &\quad +1\cdot\#\{\mathcal{C}(4,3,1)\text{ spanned by } V'\} \\
    &\quad +4\cdot\#\{\mathcal{C}(3,3,0)\sqcup\mathcal{F}(1)\text{ spanned by } V'\} \\
    &\quad +1\cdot\#\{\mathcal{C}(3,3,1)\sqcup\mathcal{F}(1)\text{ spanned by } V'\} \\
    &\quad +1\cdot\#\{\mathcal{F}(4)\text{ spanned by } V'\}.
\end{aligned}
\end{equation}

\paragraph{$p=5$.}
By Section \ref{ss:p5}, we have
\begin{equation}
\begin{aligned}
    \det L_{\omega_5}[V']&=\quad \frac{5-\sqrt{5}}{2}\cdot\#\{\mathcal{C}(4,4,0)\text{ spanned by }V'\} \\
    &\quad +\frac{5+\sqrt{5}}{2}\cdot\#\{\mathcal{C}(4,4,1)\text{ spanned by }V'\} \\
    &\quad +0\cdot\#\{\mathcal{C}(4,4,2)\text{ spanned by }V'\} \\
    &\quad +\frac{5+\sqrt{5}}{2}\cdot\#\{\mathcal{C}(4,3,0)\text{ spanned by }V'\} \\
    &\quad +\frac{5-\sqrt{5}}{2}\cdot\#\{\mathcal{C}(4,3,1)\text{ spanned by }V'\} \\
    &\quad +\frac{5+\sqrt{5}}{2}\cdot\#\{\mathcal{C}(3,3,0)\sqcup\mathcal{F}(1)\text{ spanned by }V'\} \\
    &\quad +\frac{5-\sqrt{5}}{2}\cdot\#\{\mathcal{C}(3,3,1)\sqcup\mathcal{F}(1)\text{ spanned by }V'\} \\
    &\quad +1\cdot\#\{\mathcal{F}(4)\text{ spanned by }V'\}.
\end{aligned}
\end{equation}

Remark that $\left(\mathcal{C}(3,3,g)\sqcup\mathcal{F}(1)\right)\cup\mathcal{C}(4,3,g)=\mathcal{TU}(3,3,g)\,(g\in\{0,1\})$.
Therefore
\begin{equation}
    \#\left(\mathcal{C}(3,3,g)\sqcup\mathcal{F}(1)\right)+\#\mathcal{C}(4,3,g)=\#\mathcal{TU}(3,3,g)\,(g\in\{0,1\}).
\end{equation}

By comparing these $5$ expressions,
we find out that the $5$ distinct $\omega$ values
yield a linear system of $5$ equations,
matching the $5$ types of $4$-vertex
all-regular substructures whose contributions we want to distinguish.
This system can be expressed in matrix form as follows.

\begin{equation}
\left[ \begin{array}{c}
	\det L_{\omega_2}[V']\\
	\det L_{\omega_3}[V']\\
	\det L_{\omega_4}[V']\\
	\det L_{\omega_5}[V']\\
	\det L_{\omega_6}[V']\\
\end{array} \right]=A\left[ \begin{array}{c}
	\#\mathcal{C}(4,4,0) \\
	\#\mathcal{C}(4,4,1) \\
	\#\mathcal{TU}(3,3,0) \\
	\#\mathcal{TU}(3,3,1) \\
	\#\mathcal{F}(4) \\
\end{array} \right]
\end{equation}
where
\begin{equation}
    A=\left[ \begin{matrix}
	0& 0& 4& 4&	1\\
	3& 3& 0& 3& 1\\
	0& 4& 2& 2&	1\\
	\frac{5-\sqrt{5}}{2}& \frac{5+\sqrt{5}}{2}& \frac{5+\sqrt{5}}{2}& \frac{5-\sqrt{5}}{2}& 1\\
	3& 3& 4& 1&	 1\\
\end{matrix} \right].
\end{equation}
Observe that $A$ is invertible. Therefore
\begin{equation}
\left[ \begin{array}{c}
	\#\mathcal{C}(4,4,0) \\
	\#\mathcal{C}(4,4,1) \\
	\#\mathcal{TU}(3,3,0) \\
	\#\mathcal{TU}(3,3,1) \\
	\#\mathcal{F}(4) \\
\end{array} \right]=A^{-1}
\left[ \begin{array}{c}
	\det L_{\omega_2}[V']\\
	\det L_{\omega_3}[V']\\
	\det L_{\omega_4}[V']\\
	\det L_{\omega_5}[V']\\
	\det L_{\omega_6}[V']\\
\end{array} \right],
\end{equation}
where, by calculation,
\begin{equation}
A^{-1} = \dfrac{1}{12} \left[ \begin{matrix}
-3 & -3 & -3 & 12 & -3 \\
-9 & -15 & 3 & 36 & -15 \\
-4 & -12 & 0 & 24 & -8 \\
-8 & -18 & 0 & 48 & -22 \\
60 & 120 & 0 & -288 & 120
\end{matrix} \right]+
\dfrac{\sqrt{5}}{12} \left[ \begin{matrix}
1 & 3 & -3 & 0 & -1 \\
3 & 9 & -9 & 0 & -3 \\
2 & 6 & -6 & 0 & -2 \\
4 & 12 & -12 & 0 & -4 \\
-24 & -72 & 72 & 0 & 24
\end{matrix} \right].
\end{equation}

This result demonstrates that
the numbers of these $5$ distinct types of
$4$-vertex all-regular substructures can be uniquely determined by
the $5$ values of $\det L_{\omega_i}[V']\,(i\in\{2,3,4,5,6\})$.

It is worth noting that the invertibility of matrix $A$ relies on
the distinct algebraic properties of the determinant contributions.
Notably, the row corresponding to $\omega_5$ is linearly independent
of the other rows due to its special structure,
which is essential for the full rank of the system.

This provides a linear algebra method to
precisely count these substructures within any
given $4$-vertex subset of a graph,
showcasing the versatility
of the parameterized Hermitian Laplacian framework.

\subsection{An Analog of Matrix-Tree Theorem for Parameterized Hermitian Laplacian Matrix}\label{ss:mtt}
In this section,
we present an analog of Matrix-Tree Theorem
that holds for every Hermitian Laplacian matrix
no matter the choice of the parameter $\omega$.
We start by introducing a lemma
which states that there is a bijective relationship
between spanning trees and $(n-1)$-vertex rootless forests.

\begin{lemma}\label{lem:mtt}
    If $G'(V',E')\subseteq G(V,E)$ is a regular forest substructure and
    $V'=V\setminus\{v\}$, then the substructure is
    a spanning tree of $G$ but with vertex $v$ removed,
    which is in a one-to-one correspondence with a spanning tree of $G$.
\end{lemma}

\begin{figure}[htbp]
\centering
\begin{tikzpicture}[
    vertex/.style={circle, draw, thick, minimum size=0.7cm, font=\small},
    v_prime/.style={vertex, fill=blue!20, draw=blue!50!black, very thick},
    e_prime/.style={blue!70!black, very thick, -{Stealth[length=2mm, width=1.5mm]}},
    e_prime_gray/.style={gray!70!black, very thick, -{Stealth[length=2mm, width=1.5mm]}},
    e_other/.style={gray!50, thin, -{Stealth[length=2mm, width=1.5mm]}},
    v_removed/.style={vertex, fill=gray!10, draw=gray!60, dashed},
    v_removed_added_back/.style={vertex, fill=purple!20, draw=purple!60},
    label_style/.style={font=\small}
]

\begin{scope}[xshift=0cm, yshift=0cm]
    \node[vertex] (g1) at (0,2) {1};
    \node[vertex] (g2) at (2,2) {2};
    \node[vertex] (g3) at (0,0) {3};
    \node[vertex] (g4) at (2,0) {4};
    \node[vertex] (g5) at (4,1) {5};
    \draw[e_prime_gray, shorten >=1pt, shorten <=1pt] (g1) to (g2);
    \draw[e_prime_gray, shorten >=1pt, shorten <=1pt] (g1) to (g3);
    \draw[e_prime_gray, shorten >=1pt, shorten <=1pt] (g2) to (g4);
    \draw[e_prime_gray, shorten >=1pt, shorten <=1pt] (g3) to (g4);
    \draw[e_prime_gray, shorten >=1pt, shorten <=1pt] (g4) to (g5);
    \draw[e_prime_gray, shorten >=1pt, shorten <=1pt] (g2) to (g5);
    \node[label_style, above=0.2cm of g2, align=center] {Graph $G$};
\end{scope}

\begin{scope}[xshift=6cm, yshift=0cm]
    \node[v_prime] (s1_1) at (0,2) {1};
    \node[v_prime] (s1_2) at (2,2) {2};
    \node[v_prime] (s1_3) at (0,0) {3};
    \node[v_prime] (s1_4) at (2,0) {4};
    \node[v_removed] (s1_5) at (4,1) {5};
    \draw[e_prime, shorten >=1pt, shorten <=1pt] (s1_1) to (s1_2);
    \draw[e_prime, shorten >=1pt, shorten <=1pt] (s1_1) to (s1_3);
    \draw[e_prime, shorten >=1pt, shorten <=1pt] (s1_2) to (s1_4);
    \draw[e_other, shorten >=1pt, shorten <=1pt] (s1_3) to (s1_4);
    \draw[e_other, shorten >=1pt, shorten <=1pt] (s1_4) to (s1_5);
    \draw[e_prime, shorten >=1pt, shorten <=1pt] (s1_2) to (s1_5);
    \node[label_style, above=0.2cm of s1_2, align=center] {Regular forest(tree) substructure $G'$};
    \node[below=0.1cm of s1_4, align=center, font=\tiny] {$V'=V\setminus\{5\},\,|V'|=|E'|$};
\end{scope}

\begin{scope}[xshift=3cm, yshift=-5cm]
    \node[v_prime] (t_1) at (0,2) {1};
    \node[v_prime] (t_2) at (2,2) {2};
    \node[v_prime] (t_3) at (0,0) {3};
    \node[v_prime] (t_4) at (2,0) {4};
    \node[v_removed_added_back] (t_5) at (4,1) {5};
    \draw[e_prime, shorten >=1pt, shorten <=1pt] (t_1) to (t_2);
    \draw[e_prime, shorten >=1pt, shorten <=1pt] (t_1) to (t_3);
    \draw[e_prime, shorten >=1pt, shorten <=1pt] (t_2) to (t_4);
    \draw[e_other, shorten >=1pt, shorten <=1pt] (t_3) to (t_4);
    \draw[e_other, shorten >=1pt, shorten <=1pt] (t_4) to (t_5);
    \draw[e_prime, shorten >=1pt, shorten <=1pt] (t_2) to (t_5);
    \node[label_style, above=0.2cm of t_2, align=center] {A spanning tree of $G$};
    \node[below=0.1cm of t_4, align=center, font=\tiny] {Adding vertex 5 back to $G'$ yields a spanning tree of $G$};
\end{scope}
\end{tikzpicture}
\caption{An illustration of Lemma \ref{lem:mtt}.}
\label{fig:mtt}
\end{figure}

\begin{proof}
    Obviously there can not be $3$ or more
    rootless tree components in $G'$,
    because that needs at least $2$ vertices to be dismissed,
    but $v$ is the only dismissed vertex.

    Suppose adding $v$ to the vertex subset $V'$.
    Doing so turns one rootless tree (or both two rootless trees)
    in $G'$ into one normal tree.
    No matter which case it is, this turns $G'$
    into a spanning tree of $G$.
\end{proof}

By this fact we get a sufficient condition of the modulus
of any cofactor of a Hermitian Laplacian matrix
of an oriented graph $G$ being equal to the number of spanning trees
of its underlying graph.

\begin{theorem}\label{thm:mtt}
    Suppose $G(V,E)$ is an oriented graph
    whose Hermitian Laplacian matrix is $L_\omega$.
    If under the parameter $\omega$,
    among all $(n-1)$-vertex substructure of $G$,
    there are no unicyclic graph components
    that have a non-zero Hermitian Laplacian determinant,
    then any cofactor of $L_\omega$
    is equal to the number of spanning trees
    of the underlying graph of $G$.
\end{theorem}

\begin{proof}
    Arbitrarily pick $v\in V$ and let $V'=V\setminus\{v\}$.
    By Theorem \ref{thm:lapl-decomp}:
    \begin{equation}
        \det L_\omega[V']=\sum_{\substack{E':(V',E')\text{ all-regular}\\|E'|=|V'|=|V|-1}}\det L_\omega(G'(V',E')).
    \end{equation}

    By the given condition,
    each all-regular substructure
    that has non-zero Hermitian Laplacian determinant
    must be a regular forest,
    therefore
    \begin{equation}
    \begin{aligned}
        \det L_\omega[V']&=\sum_{E':(V',E')\text{ is a }(|V|-1)\text{-vertex regular forest}}\underbrace{\det L_\omega(G'(V',E'))}_{=1} \\
        &=\text{number of }(n-1)\text{-vertex regular forest spanned by }G \\
        &=\text{number of spanning trees of }G.
    \end{aligned}
    \end{equation}

    This completes the proof.
\end{proof}

\begin{remark}
    If we take $\omega=1$, then every cycle has a
    Hermitian Laplacian determinant of $0$.
    In this case,
    no unicyclic graph has non-zero Hermitian Laplacian determinant,
    therefore any cofactor of the Hermitian Laplacian matrix
    will be equal to the number of the spanning trees.
    This is the well-known Kirchoff's Theorem [\citenum{kirchoff}].

    Theorem \ref{thm:mtt} suggests that
    if we can select certain parameter $\omega$
    and set constraints on the graphs
    such that every cycle in such graphs
    has a determinant of $0$ under the parameter,
    then any cofactor of this Hermitian Laplacian matrix
    can be equal to the number of its spanning trees.
    For example, Xiong et al. provided another analog
    of Matrix-Tree Theorem ([\citenum{sk-2}], Theorem 3.7)
    for the special case of $\omega=\ee^{\im 2\pi/6}$,
    where they introduced a special graph type
    called S P graph for the theorem to hold under.
\end{remark}

\begin{proposition}
    Let $\omega$ be a complex number with modulus $1$.
    Suppose $G(V,E)$ is an oriented graph
    such that every unicyclic subgraph of $G$
    has zero Hermitian Laplacian determinant
    under parameter $\omega$.
    Denote $L_\omega$ to be the Hermitian Laplacian matrix of $G$.
    Then any cofactor of $L_\omega$
    is equal to the number of spanning trees
    of the underlying graph of $G$.
\end{proposition}

This section illustrates another application of
the parameterized Hermitian Laplacian framework.

\section{Conclusion}\label{sec:conclusion}

In this paper we presented an analysis
of oriented graph substructures
using a parameterized Hermitian Laplacian matrix $L_\omega$.
The approach is based on selecting parameters
by leveraging the methodology of Galois theory,
particularly focusing on primitive roots of unity
and their Galois conjugates.
We established through the
Hermitian Laplacian Principal Minor Decomposition Theorem
that principal minors of $L_\omega$ can be expressed
as sums over the determinants of all-regular substructures,
specifically rootless forests and unicyclic graphs.

A central outcome of this work is a general method
for enumerating non-$p$-vanishing unicyclic graph components
within a given all-regular substructure, for any odd prime $p$.
This method arises from our key theoretical finding:
the product of Hermitian Laplacian determinants,
parameterized over a full set of Galois conjugates
of a primitive $p$th root of unity,
equals an integer power of $p$.

As a concrete and illustrative example,
we detailed the $p=5$ case,
where the determinants lie in the field $\Q(\sqrt{5})$,
and showed how their product and ratio
can directly yield the counts of two different types
of unicyclic graph components.
The utility of the variable parameter was also shown
in developing a method to count
distinct $4$-vertex all-regular substructures
using a linear system of equations
derived from different roots of unity.

Several avenues for future investigation remain.
Exploring more roots of unity with special properties
or parameters within different algebraic number fields
might allow for the distinction of more complex substructures.
In addition, extending this parameterized approach
to other graph matrices,
such as Hermitian quasi-Laplacian matrices,
could also be considered.
Furthermore, the development of algorithms
based on the theoretical findings of this paper
could provide practical tools for graph analysis.

In summary, by examining the connections between
the algebraic properties of
the parameterized Hermitian Laplacian matrix
and the topology of oriented graph substructures,
this research may provide a new tool
for algebraic graph theory
and the structural investigation of graphs and networks.

\section*{Acknowledgements}
The author wishes to express gratitude to his advisor,
Prof. Gehao Wang,
for his valuable writing advice,
and to the anonymous reviewers
for their constructive comments
on this article.

\end{document}